%% file: paper_sdf-usmnet.tex
\documentclass[a4paper]{article}
\usepackage{authblk}
\input{packages.tex}
\input{macros.tex}

\title{\papertitle}
\date{}

\author[1,2]{Linying Zhang}
\author[1]{Stefano Pagani}
\author[2]{Jun Zhang}
\author[1]{Francesco Regazzoni
\footnote{\textit{Corresponding author} (\texttt{francesco.regazzoni@polimi.it})}}

\affil[1]{\small{MOX, Laboratory of Modeling and Scientific Computing, Dipartimento di Matematica, Politecnico di Milano, Milano, Italy}}
      
\affil[2]{\small{School of Aeronautic Science and Engineering, Beihang University, Beijing, China}}

\begin{document}

\maketitle

\begin{abstract}
    \input{parts_abstract}
\end{abstract}


\input{parts_intro}
\input{parts_method}
\input{parts_data_preparation}
\input{parts_results}

\input{parts_discussion}

\input{parts_appendix}

\section*{Acknowledgement}
\input{parts_acknowledgements}

\bibliographystyle{plain}
\bibliography{ref}

\end{document}

%% file: packages.tex
\usepackage{graphicx} 
\usepackage{amsmath}
\usepackage{amssymb}
\usepackage{mathrsfs}
\usepackage{subcaption}
\usepackage{caption}
\usepackage{tabularx}
\usepackage{float}
\usepackage{newfloat}
\usepackage{algpseudocode}
\usepackage{algorithm}
\usepackage{todonotes}
\usepackage{enumitem}

%% file: macros.tex
\newcommand{\papertitle}{Shape-informed surrogate models based on signed distance function domain encoding}

\newcommand{\keywordOne  }{Scientific Machine Learning}
\newcommand{\keywordTwo  }{Surrogate Modeling}
\newcommand{\keywordThree}{Shape Models}
\newcommand{\keywordFour }{Neural Networks}
\newcommand{\keywordFive }{Signed Distance Functions}

\DeclareMathOperator*{\argmax}{arg\,max}
\DeclareMathOperator*{\argmin}{arg\,min}
\DeclareFloatingEnvironment[
    fileext=lob,
    listname={List of Boxes},
    name=Box,
    placement=tbhp,
    within=section, 
]{mybox}

\newcommand{\nnSDF}{\mathcal{NN}_\text{SDF}}
\newcommand{\nnPHY}{\mathcal{NN}_\text{Phys}}

\newcommand{\paramsSDF}{\boldsymbol{\theta}_\text{SDF}}
\newcommand{\paramsPHY}{\boldsymbol{\theta}_\text{Phys}}
\newcommand{\paramsSDFhat}{\hat{\boldsymbol{\theta}}_\text{SDF}}
\newcommand{\paramsPHYhat}{\hat{\boldsymbol{\theta}}_\text{Phys}}

%% file: parts_abstract.tex
We propose a non-intrusive method to build surrogate models that approximate the solution of parameterized partial differential equations (PDEs), capable of taking into account the dependence of the solution on the shape of the computational domain. Our approach is based on the combination of two neural networks (NNs). The first NN, conditioned on a latent code, provides an implicit representation of geometry variability through signed distance functions. This automated shape encoding technique generates compact, low-dimensional representations of geometries within a latent space, without requiring the explicit construction of an encoder. The second NN  reconstructs the output physical fields independently for each spatial point, thus avoiding the computational burden typically associated with high-dimensional discretizations like computational meshes. Furthermore, we show that accuracy in geometrical characterization can be further enhanced by employing Fourier feature mapping as input feature of the NN. The meshless nature of the proposed method, combined with the dimensionality reduction achieved through automatic feature extraction in latent space, makes it highly flexible and computationally efficient. This strategy eliminates the need for manual intervention in extracting geometric parameters, and can even be applied in cases where geometries undergo changes in their topology. Numerical tests in the field of fluid dynamics and solid mechanics demonstrate the effectiveness of the proposed method in accurately predict the solution of PDEs in domains of arbitrary shape. Remarkably, the results show that it achieves accuracy comparable to the best-case scenarios where an explicit parametrization of the computational domain is available.

%% file: parts_intro.tex
\section{Introduction}\label{sec:intro}

Reduced order models (ROMs) provide efficient and reliable approximations of high-fidelity models, known as full-order models (FOMs), which are often based on methods like the Finite Element Method (FEM) or the Finite Volume Method (FVM). ROMs are designed to offer a significantly reduced computational cost compared to the corresponding FOMs.
ROMs can be classified as either intrusive or non-intrusive. Intrusive ROMs require explicit knowledge of the FOM, such as in projection-based methods \cite{sirovich1987turbulence,quarteroni2015reduced,hesthaven2016certified}. Non-intrusive ROMs, also referred to as surrogate models or emulators, are constructed from a set of input-output pairs generated by the FOM, either without knowledge of the underlying equations or with only partial understanding of the physics involved \cite{peherstorfer2015dynamic,regazzoni2019modellearning,guo2019data,vlachas2022multiscale}.

In recent years, there has been a notable rise in methods for constructing surrogate models using machine learning tools such as neural networks (NNs), Gaussian processes, and support vector machines \cite{brunton2020machine,lee2020model,longobardi2020predicting,maulik2021reduced,oommen2022learning,regazzoni2024ldnets}. 
Despite this progress, most existing methods, whether intrusive or non-intrusive, are tailored to a specific geometry and require re-training when the domain shape changes.
This limitation is particularly challenging in various applications, including the biomedical field, where the solutions of mathematical models are highly dependent on the shape of the computational domain, which can vary significantly from patient to patient.

Efforts to handle changes in the shape of domains flexibly have leveraged geometrical shape models, which represent different geometries by varying a limited number of parameters. Some of the most commonly used techniques include:
free-form deformations~\cite{sederberg1986free,lamousin1994nurbs}, where points in the domain are deformed via interpolation of displacements on a lattice grid;
radial-basis functions~\cite{buhmann2000radial}, which define a map parameterized by the displacement of a few control points;
statistical shape models~\cite{heimann2009statistical}, such as those based on principal component analysis (PCA)~\cite{jolliffe2016principal}, which describe the variation of geometries as a weighted sum of modes of variation with respect to an average shape, possibly combined with Large-Deformation Diffeomorphic Metric Mapping~\cite{ardekani2009computational};
isogeometric analysis (IGA) which, similarly, handles changes in geometry without the need to regenerate a new mesh, by using non-uniform rational B-splines (NURBS)~\cite{hughes2005isogeometric}.

In the context of model order reduction, shape models have been employed to construct geometric parameterizations that can be explicitly handled by the ROMs in the same way as physical parameters are treated \cite{lassila2010parametric,manzoni2012model,sangalli2009case}. 
However, these methods often depend on deforming a predefined mesh, which restricts the range of geometries that can be effectively handled to guarantee sufficient mesh quality. 
Additionally, in practical scenarios, statistical shape models like PCA may need a large number of modes to accurately capture geometric variability \cite{romdhani1999multi}. 
Moreover, even when shape models significantly reduce the dimensionality of geometric variability, they often fall short in effectively extracting useful features for predicting the solution of differential models defined on the geometry considered (see, e.g., \cite{rodero2021linking}).

This complexity has motivated researchers to shift to geometrically informed deep learning models for reduced-order modeling in variable shape.
One architecture designed to manage changes in mesh topology is the Graph Neural Network (GNN) \cite{wu2020comprehensive}, which has been employed to train surrogate models that can generalize to unseen geometries \cite{pegolotti2024learning}.
PointNets \cite{qi2017pointnet}, namely NNs that can handle point-cloud representations of geometries, thanks on invariant operations with respect to point ordering and on local-global aggregation, provide another tool to inform surrogate models with geometrical features \cite{oldenburg2022geometry}.
Convolution-deconvolution operations have been applied to perform an image-to-image regression \cite{guo2016convolutional}, which however suffers from inaccuracy near the boundary due to the lack of very fine pixels in the geometry images.
A recently proposed method addressing limited resolution challenges is USM-Net \cite{regazzoni2022universal}, namely a surrogate modeling technique which applies to differential problems whose solution depends on both physical and geometrical parameters. 
A USM-Net is based on two main components.
The first component is a shape-encoder, which is a map that associates to each geometry a (low-dimensional) vector of scalars, called \textit{shape codes} or \textit{geometrical features}, encoding the main properties of the considered geometry.
The second component is a coordinate-based NN (that is a NN taking space coordinates as inputs, also known as neural field \cite{xie2022neural}), conditioned on both the physical parameters and the geometrical features, that predicts the output of the PDE in each point of the domain.
Thanks to their mesh-less nature, USM-Nets allow to flexibly handle changes of shape of the domain, while relying on lightweight architectures that yield very good generalization properties (see \cite{regazzoni2024ldnets} for an investigation on the effectiveness of coordinate-based NNs to approximate the solution of PDEs).
Optionally, USM-Nets can be supplemented with a universal coordinate system, consisting in a map of all the different geometries onto a reference one, thus enhancing the surrogate model accuracy \cite{regazzoni2022universal}.

In \cite{regazzoni2022universal}, two shape-encoding techniques were considered. 
The first technique, applicable when an exact parametrization of the geometries is available, simply consists in employing the parameters used to construct the geometry as shape codes.
This approach, despite being highly appealing in contexts like industrial design, is not practical in fields such as geosciences and biomedicine, that is when the domain is not human-designed and an explicit parametrization is unavailable.
The second shape-encoding technique, instead, defines the shape codes as the collection of landmarks, corresponding to coordinates of key points of the domain.
This approach too features some limitations: it requires a manual intervention for tagging these points and a case-by-case definition, which is unlikely to identify the minimum set of points most informative for the USM-Net about geometric variability. As a result, it is hardly generalizable.

To address these challenges, the goal of this work is to develop an enhanced version of USM-Nets, based on a shape encoding technique that is fully automatic and applicable to a wide range of cases. 
The proposed method, named SDF-USM-Nets, leverages DeepSDF \cite{park2019deepsdf}, a NN-based architecture that is able to learn nonlinear shape models based on implicit representations. 
Specifically, DeepSDFs are trained to approximate the signed distance function (SDF) of the geometries belonging to the training set through a coordinate-based NN, conditioned on a latent code.
By combining these two approaches, SDF-USM-Nets are equipped with these key innovations:
\begin{itemize}
    \item Shape encoding is directly provided by the latent code of a geometrical DeepSDF model, removing the need for user-based geometrical feature extraction.
    \item The SDF field, employed to define the shape codes, also offers a convenient way of generalizing the concept of the universal coordinate system introduced in \cite{regazzoni2022universal}, by providing the USM-Net with an additional feature that improves, as testified by the numerical tests, its generalization properties across different shapes.
    \item Remarkably, this shape encoding technique is compatible with large deformations and topological changes, thus offering great flexibility.
\end{itemize}
The paper is organized as follows. 
First, in Sec.~\ref{sec:methods}, we introduce the notation considered in this paper and we present the proposed method.
Then, in Sec.~\ref{sec:data}, we detail the test cases considered in this work, and in Sec.~\ref{sec:results} we provide the associated numerical results.
Finally, in Sec.~\ref{sec:discussion} we discuss them and provide some final remarks.

%% file: parts_method.tex
\section{Methods} \label{sec:methods}

In this Section, we first define the problem of interest. Then, we introduce the proposed SDF-USM-Net architecture and the training strategy.

\subsection{Problem definition}
We define \(\mathbf{u}:\Omega \rightarrow \mathbb{R}^{d_\mathbf{u}}\) as the quantity of interest (QoI), which may correspond, to provide a couple of practical examples, to the velocity field in a vessel or the displacement field of soft material. This QoI is defined in the spatial domain \(\Omega \subset \mathbb{R}^d\), where typically \(d=2,3\). 
As a paradigmatic (albeit not unique) case, the QoI might be governed by a PDE in the following form:
\begin{equation}\label{eq1}
\begin{aligned}
    \mathcal{N}(\mathbf{u})&=0,\ &&\text{in } \Omega, \\
    \mathcal{B}(\mathbf{u})&=0,\ &&\text{on } \partial \Omega,  
\end{aligned}
\end{equation}
where \(\mathcal{N}\) is a generic differential operator and \(\mathcal{B}\) is the operator associated with boundary conditions, defined on the boundary \(\partial \Omega\) of the domain. 

In numerous practical problems, the solution field \(\mathbf{u}(\mathbf{x})\) is intrinsically dependent on the spatial domain \(\Omega\) itself. 
We will employ the notation \(\mathbf{u}(\mathbf{x}; \Omega)\) to stress this aspect.
With a well-posed differential problem \eqref{eq1}, there exists a unique solution \(\mathbf{u}(\mathbf{x}; \Omega)\), which can be numerically approximated through a full-order model (FOM). Our focus lies in exploring the manifold of fields \(\mathbf{u}(\cdot;\cdot)\) generated upon variations of the domain \(\Omega\), and learning a solution map \(\Omega \to \mathbf{u}(\cdot; \Omega)\) between the variable domain \(\Omega\) and the QoI \(\mathbf{u}(\mathbf{x}; \Omega)\) from a finite collection of \(\{\Omega, \mathbf{u}(\cdot; \Omega)\}\) pairs, sampled from a predetermined subset of interest for the practical problems at hand, denoted by \(\mathcal{G}\):
\[\mathcal{G} \subset \{\Omega \subset \mathbb{R}^d,\ \mathrm{open\ and\ bounded}\}.\]

\subsection{SDF-USM-Net}
\begin{figure}
    \centering
    \includegraphics[width=0.7\linewidth]{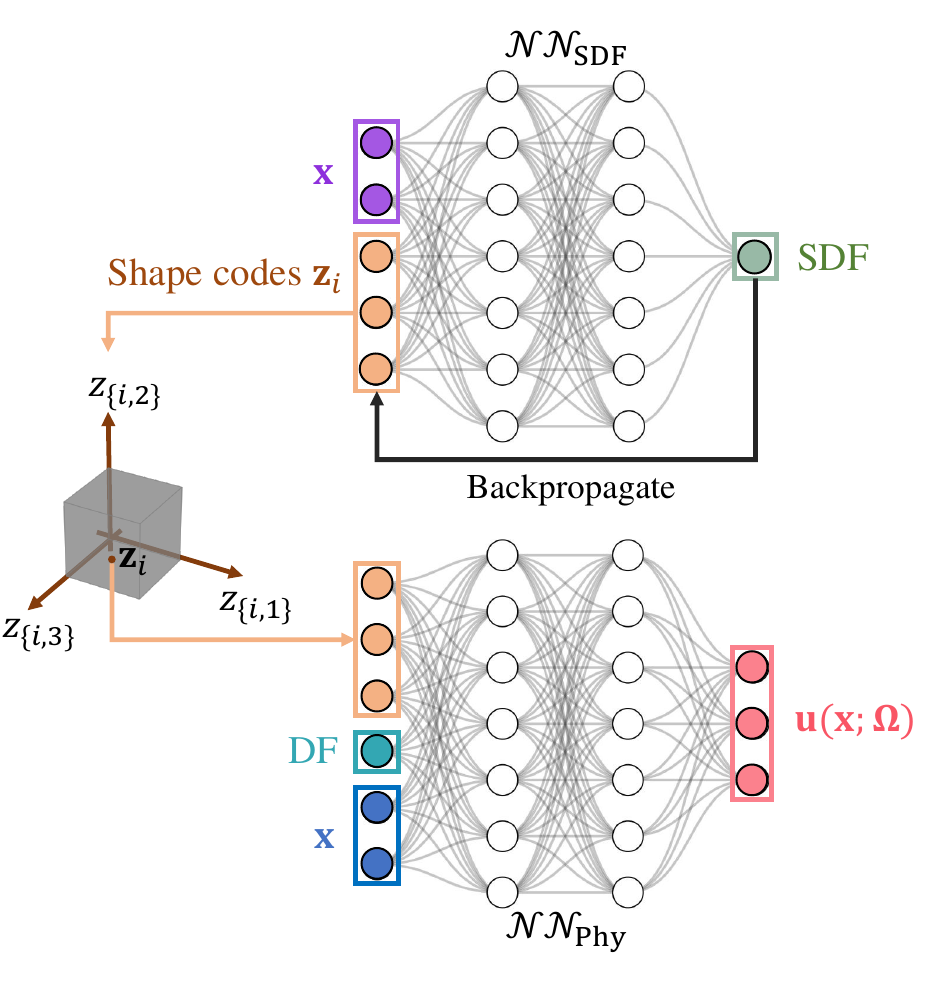}
    \caption{SDF-USM-Net architecture. The network \(\nnSDF\) receives the input spatial coordinates \(\mathbf{x}\) and the shape code \(\mathbf{z}_i\) and outputs the SDF value for the corresponding point. Within an auto-decoder framework, it learns the shape code distribution in a latent space by reconstructing the continuous SDF in the spatial domain. After the training of \(\nnSDF\), the network \(\nnPHY\) is trained in a supervised manner, taking as input the shape code \(\mathbf{z}_i\), spatial coordinates \(\mathbf{x}\), and a distance function \(DF(\mathbf{x})\) value and returning the approximation of QoI \(\mathbf{u}(\mathbf{x};\Omega)\).} 
    \label{struc}
\end{figure}

The architecture of our proposed SDF-USM-Net, shown in Fig. \ref{struc}, consists of two sub-networks, \(\nnSDF(\mathbf{x},\mathbf{z}_i;\paramsSDF)\) and \(\nnPHY(\mathbf{x},\mathbf{z}_i,DF(\mathbf{x});\paramsPHY)\), which are two fully connected neural networks (FCNN) with trainable parameters \(\paramsSDF\) and \(\paramsPHY\), respectively, and a trainable shape code vector \(\mathbf{z}_i \in \mathbb{R}^k\) for each geometry. 
Here, \(DF(\mathbf{x})\) denotes the distance from a point \(\mathbf{x}\) to the nearest surface. 
An SDF-USM-Net defines a map from the domain \(\Omega \in \mathcal{G}\) and the spatial coordinates input \(\mathbf{x}\in \Omega\) to the QoI \(\mathbf{u}\) associated with the input query point \(\mathbf{x}\) and geometry \(\Omega\). 
The first NN, namely \(\nnSDF\), allows to associate to any shape \(\Omega_i \in \mathcal{G}\) a finite-dimensional shape code vector \(\mathbf{z}_i\) to give a compact description of \(\Omega_i\). The second NN, \(\nnPHY\), associates to the shape code \(\mathbf{z}\) the solution in terms of physical field \(\mathbf{u}(\mathbf{x}; \Omega)\) corresponding to the input query point \(\mathbf{x}\).
The two NNs are trained sequentially.

\subsubsection{Geometry feature extraction}\label{sec:deepsdf}
\begin{figure}
    \centering
    \includegraphics[width=0.9\linewidth]{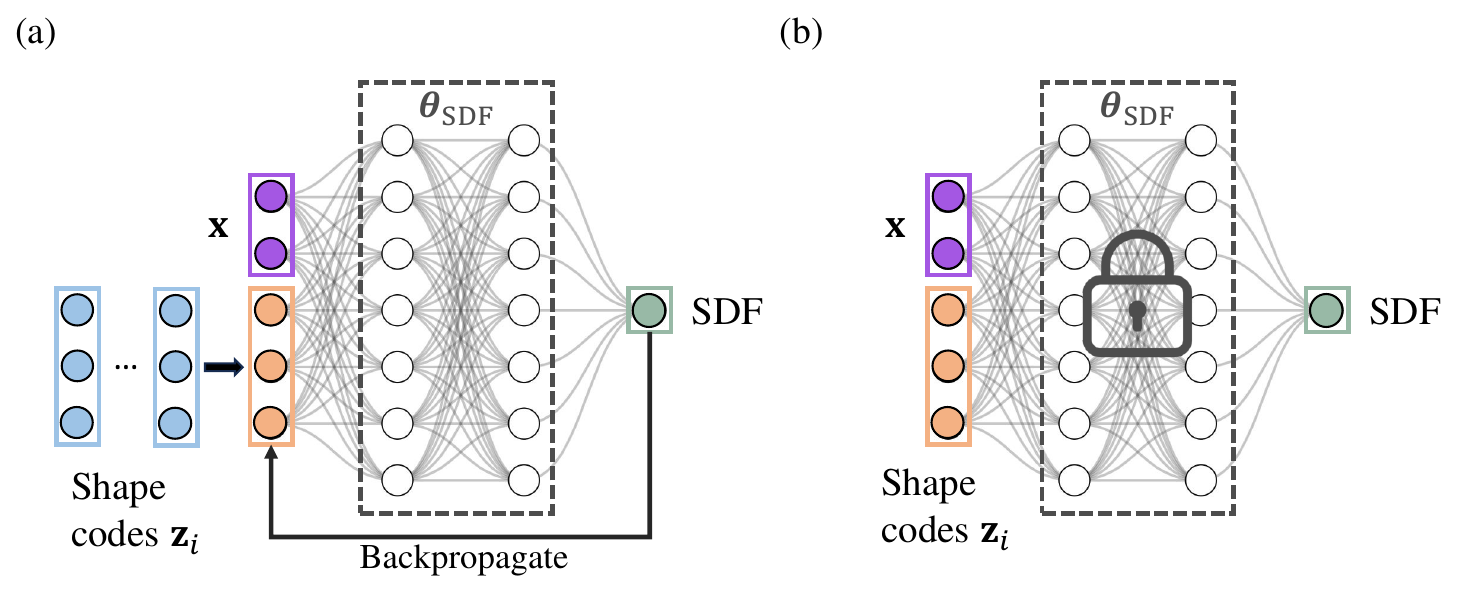}
    \caption{
    Working principle of \(\nnSDF\). (a) During training, each sample geometry is associated with a randomly initialized shape code. The auto-decoder receives shape codes and coordinates as input. The shape codes are optimized along with the decoder trainable parameters \(\paramsSDF\) through standard back-propagation. (b) During inference, the trainable parameters \(\paramsSDF\) are fixed, and an optimal shape code \(\mathbf{z}_i\) is estimated for any geometry by minimizing the discrepancy between predicted and expected SDF values.}
    \label{fig:deepsdf}
\end{figure}

The construction of \(\nnSDF\) is based on the concept of DeepSDF \cite{park2019deepsdf}, an implicit function-based shape representation method, that we use here to embed the space of shapes into a low-dimensional latent space (see Figure \ref{fig:deepsdf} for a graphical explanation of the model).

In this context, signed distance functions (SDFs) are utilized to describe shapes in the Euclidean space, so that the zero level-set of the function implicitly defines the domain boundary. The SDF is defined as follows
\begin{equation} \label{eq2}
    SDF(\mathbf{x})=
    \begin{cases}
        d(\mathbf{x}, \partial \Omega) & \text{if $\mathbf{x} \in \Omega$,} \\
        - d(\mathbf{x}, \partial \Omega) & \text{if $\mathbf{x} \notin \Omega$,} \\
    \end{cases}
\end{equation}
where \(d\) denotes the minimum distance of a point $\mathbf{x}$ from the boundary $\partial \Omega$. 
A NN, denoted as \(\nnSDF\), is trained on a set of \(N\) shapes \(\{G_i\}_{i=1}^N\). Each shape \(G_i\) is represented as a set of data pairs \(G_i:=\{(\mathbf{x}_i^j, s_i^j):s_i^j=SDF_i(\mathbf{x}_i^j)\}_{j=1}^{K_i}\), where \(\mathbf{x}_i^j\) is sampled around the boundary of the shape \(G_i\), and \(s_i^j\) is the precomputed SDF value ($SDF_i$ denotes the SDF associated with the domain $\Omega_i$). 

With an FCNN architecture, \(\nnSDF\) receives as input the space coordinate \(\mathbf{x}\) and a trainable vector \(\mathbf{z}_i\in \mathbb{R}^k\) which will be later interpreted as the shape code that represents the target shape in the \(k\)-dimensional latent space. The FCNN output is an estimate of the SDF value associated with the input coordinate \(\mathbf{x}\). The dimension of the latent space \(k\) is set by the user and should be regarded as a hyperparameter. The trainable parameters \(\paramsSDF\) of FCNN and shape codes \(\{\mathbf{z}_i\}_{i=1}^N\) are learned from the training set \(G^{\mathrm{train}}:=\{G_i\}_{i=1}^N\) to make the output of the NN a good approximator of the given SDF of the target domain \(\Omega\):
\begin{equation}
    SDF_i(\mathbf{x}) \approx \nnSDF(\mathbf{x}, \mathbf{z}_i; \paramsSDF).
\end{equation}
More precisely, the trainable parameters \(\paramsSDF\) are learned simultaneously to the shape codes \(\{\mathbf{z}_i\}_{i=1}^N\) associated with training shapes \(G_i\) within an \textit{auto-decoder} framework.
Hence, by following \cite{park2019deepsdf}, we consider a zero-mean multivariate-Gaussian \(\mathcal{N}(0, \frac{1}{k} I)\) as the prior distribution \(p(\mathbf{z}_i)\) of shape codes, where \(k\) is the dimension of latent space, and we compute the parameters \(\paramsSDF\) and \(\mathbf{z}_i\) that maximize the posterior across all shapes:
\begin{equation} \label{eq:argmin}
\begin{split}
    \paramsSDFhat & =\argmax_{\paramsSDF} \sum_{i = 1}^N \max_{\mathbf{z}_i} \log( p_{\paramsSDF} (\mathbf{z}_i|G_i) ) \\
   & =\argmax_{\paramsSDF} \sum_{i = 1}^N \max_{\mathbf{z}_i} \left[\log( p_{\paramsSDF} (G_i|\mathbf{z}_i) )+\log( p(\mathbf{z}_i) )\right]
   \\
   & =\argmax_{\paramsSDF} \sum_{i = 1}^N \max_{\mathbf{z}_i} \left[\sum_{j=1}^{K_i}\log( p_{\paramsSDF} (s_i^j|\mathbf{z}_i) )+\log( p(\mathbf{z}_i) )\right].
\end{split}
\end{equation}
We assume that the likelihood takes the form
\begin{equation}
    p_{\paramsSDF} (s_i^j|\mathbf{z}_i) 
    = \exp\{-\mathcal{L}(
        s_i^j, \,
        \nnSDF(\mathbf{x}_i^j, \mathbf{z}_i; \paramsSDF)
    )\},
\end{equation}
where \(\mathcal{L}\) is a loss function penalizing the difference between the predicted SDF and the reference one.
Specifically, we consider two options.
The first one consists in the $L_1$ loss function 
$\mathcal{L}(s, \tilde{s})=|s - \tilde{s}|$, where $s$ denotes the reference SDF and $\tilde{s}$ the approximated one.
The second option consists in the smoothly clamped $L_1$ loss function:
\begin{equation} \label{eqn:clampL1}
    \mathcal{L}(s, \tilde{s})=|\mathrm{clamp}(s,\beta)-\mathrm{clamp}(\tilde{s},\beta)|,
\end{equation}
where, \(\mathrm{clamp}_\beta(s):=\beta\ \mathrm{tanh}(s/\beta)\) is a function to smoothly truncate the SDF values at a controlled distance, aimed at focusing the training on the details of the boundary.
Notice that, in the original DeepSFD work \cite{park2019deepsdf}, a sharply truncated $L_1$ loss function is employed, corresponding to \eqref{eqn:clampL1} where the above definition is replaced by \(\mathrm{clamp}_\beta(s):=\min(\beta, \max(-\beta, s))\).
Our tests suggested that the use of a smoothed clamping is helpful is escaping local minima, due to the lack a flat profile when the absolute value of the argument is larger than $\beta$.
Furthermore, the use of the truncation function itself may be convenient or not depending upon the specific application. 
As a rule of thumb, when the geometric dataset under consideration manifests numerous variations in detail, opting to utilize the truncation function may be convenient. Conversely, in cases when the geometric dataset features alterations in the topology of shapes, abstaining from the truncation function likely renders the shape code inference more stable upon topological changes.
In Sec.~\ref{sec:results} we numerically investigate this aspect.

In conclusion, the final cost function that is minimized with respect to the network parameters \(\paramsSDF\) and shape codes \(\mathbf{z}_i\) reads:
\begin{equation} \label{eqn:SDF_training}
\paramsSDFhat, \{\hat{\mathbf{z}}_i\}^N_{i=1} = 
    \argmin_{\paramsSDF,\{\mathbf{z}_i\}^N_{i=1}} \sum_{i=1}^{N} 
    \left(\sum_{j=1}^{K_i} \mathcal{L}(
        s_i^j, \,
        \nnSDF(\mathbf{x}_i^j, \mathbf{z}_i; \paramsSDF)
    )+\frac{1}{\sigma^2}||\mathbf{z}_i||^2_2 \right),
\end{equation}
where the regularization parameter \(\sigma\) balances in practice the SDF reconstruction and the sparsity of codes.


Once the parameters $\paramsSDFhat$ of $\nnSDF$ are trained, the auto-decoder can be used to get the shape code of a new geometry by minimizing the discrepancy between the SDF values obtained from $\nnSDF$ and observed SDF samples.
More precisely, suppose that a new geometry in the test dataset $G^{\mathrm{test}}$ is given in the form of 
\begin{equation*}
G_i:=\{(\mathbf{x}_i^j, s_i^j):s_i^j=SDF_i(\mathbf{x}_i^j)\}_{j=1}^{K_i}, \quad i > N.
\end{equation*}
Then, the associated shape code \(\mathbf{z}_i\) can be obtained through:
\begin{equation} \label{eqn:shape_code_inference}
    \hat{\mathbf{z}}_i=\argmin_{\mathbf{z}} \sum_{j=1}^{K_i} \mathcal{L}(
        s_i^j, \,    
        \nnSDF(\mathbf{x}_i^j, \mathbf{z}; \paramsSDFhat)
    ) +
    \frac{1}{\sigma^2}||\mathbf{z}||^2_2.
\end{equation}
In this stage, the shape code \(\mathbf{z}_i\) is randomly initialized from the empirical statistical multivariate-Gaussian \(\mathcal{N}(\hat{\boldsymbol{\mu}}, \hat{\Sigma})\), where $\hat{\boldsymbol{\mu}}$ and $\hat{\Sigma}$ are the sample mean and covariance of the population of the $N$ learned shape codes $\hat{\mathbf{z}}_i$ (in the training set \(G^{\mathrm{train}}\)),
to facilitate a swift convergence towards an optimal \(\mathbf{z}_i\).

\subsubsection{Physical field prediction}

After the \(\nnSDF\) training stage, consisting of learning the distribution of shape codes in the latent space through a set of SDFs, a corresponding shape code vector \(\mathbf{z}_i\) can be obtained for each geometry \(\Omega_i \subset \mathcal{G}\) (see Fig. \ref{struc}). We then leverage this association to train a second NN, denoted as \(\nnPHY\), to map each geometry to the QoI \(\mathbf{u}(\mathbf{x}; \Omega)\).

\(\nnPHY\) is trained based on a set of $M$ samples \(\{U_i\}_{i=1}^M\), in which $M$ can in principle be different from the number of training samples $N$ used in \(\nnSDF\).
Each physical field in the training set is represented by $P_i$ data points in a set \(U_i := \{(\mathbf{x}_i^j, \overline{s}_i^j, \mathbf{z}_i), \mathbf{u}(\mathbf{x}_i^j;\Omega_i)\}_{j=1}^{P_i}\), where \(\mathbf{u}(\mathbf{x}_i^j;\Omega_i)\) is the FOM solution. 
Different from the data pairs considered to train \(\nnSDF\), \(\mathbf{x}_i^j\) is sampled inside the computational domain.
By \(\overline{s}_i^j\) we denote the distance of each point from the domain physical boundary, with the aim to provide the NN with an additional feature to improve generalization among different geometries.
Hence, depending on the problem at hand, different definitions for \(\overline{s}_i^j\) can be considered, in order to better reflect the underlying physics.
Practical examples will be provided in Sec.~\ref{sec:results}.
In general, we define \(\overline{s}_i^j = {DF}_i(\mathbf{x}_i^j)\), namely the distance function ${DF}_i(\mathbf{x}) = d(\mathbf{x}, \Gamma_i)$, where $\Gamma_i$ is a suitable subset of $\partial\Omega_i$.
In the case when $\Gamma_i = \partial\Omega_i$, we can set ${DF}_i \equiv {SDF}_i$, as only points internal to $\Omega_i$ are considered.
Each data pair $\{\mathbf{x}_i^j, \overline{s}^j_i\}$ is stacked with the shape codes \(\mathbf{z}_i\). If the considered geometry is among those used to train $\nnSDF$ $(i\leq N)$, then the shape codes obtained by \eqref{eqn:SDF_training} are employed. Otherwise $(i>N)$, shape codes obtained through \eqref{eqn:shape_code_inference} are used.

Based on an FCNN architecture, \(\nnPHY\) is trained in a supervised manner based on the set \(\{U_i\}_{i=1}^M\), in order to approximate the FOM solutions as:
\begin{equation}
\mathbf{u}(\mathbf{x};\Omega_i) \approx 
\nnPHY(\mathbf{x},\mathbf{z}_i, \overline{s}_i^j(\mathbf{x}); \paramsPHY),
\end{equation}
More precisely, the trainable parameters $\paramsPHY$ are obtained by minimizing the following mean-square error (MSE) cost function:
\begin{equation}\label{eq:phy_training}
\paramsPHYhat = \argmin_{\paramsPHY} 
\sum_{i=1}^{M} \sum_{j=1}^{P_i} 
| 
\mathbf{u}(\mathbf{x}_i^j; \Omega_i) -  
\nnPHY(\mathbf{x}_i^j,{\mathbf{z}}_i, \overline{s}_i^j; \paramsPHY)
|_2^2.
\end{equation}
After obtaining the parameters \(\paramsSDFhat\) and \(\paramsPHYhat\) from the training of $\nnSDF$ and $\nnPHY$, respectively, in the inference stage, the solution \(\mathbf{u}(\mathbf{x}; \Omega_i)\) for a new geometry can be approximated through the following two steps: 
\begin{enumerate}[label=(\roman*)]
    \item compute the shape code \(\hat{\mathbf{z}}_i\), where \(i > M\), by minimizing Eq.~\eqref{eqn:shape_code_inference}; 
    \item evaluate the solution \(\mathbf{u}(\mathbf{x}; \Omega_i)\) through a single feedforward pass in $\nnPHY$ for a given query point \(\mathbf{x}\) together with the shape code \(\hat{\mathbf{z}}_i\) and the associated value $\overline{s}$.
\end{enumerate}

\subsubsection{Centralization for discrete shape learning} \label{sec: central}
To facilitate the learning of shapes, a centralization step may be implemented (this is exemplified in Test Case 2). 
In this case, the input coordinates \(\mathbf{x}\) to \(\nnSDF\) are centralized using the equation:
\[
\hat{\mathbf{x}} = \mathbf{x} - \mathbf{x_0},
\]
where \(\mathbf{x_0}\) represents the centroid of the shape. Subsequently, the centroid \(\mathbf{x_0}\) is concatenated to the shape code, and thus provided as an additional input to \(\nnPHY\). 

\subsubsection{Normalization and input mapping} \label{sec:mapping}
To facilitate the training of the NNs, we normalize the input coordinates \(\mathbf{x}\) of both \(\nnSDF\) and \(\nnPHY\) and the output physical fields \(\mathbf{u}\) to the interval [-1, 1]. Specifically, by considering the interval \([\mathbf{\alpha}_{\min}, \mathbf{\alpha}_{\max}]\), we preprocess each scalar variable \(\alpha\) through the following normalization layer  
\begin{equation}
    \alpha \mapsto 2\frac{\alpha-\alpha_{\min}}{\alpha_{\max}-\alpha_{\min}} - 1.
\end{equation}
Moreover, we optionally use Fourier feature mapping \cite{tancik2020fourier} \(\gamma\) to preprocess input coordinates before passing them through the FCNN. The function \(\gamma\) maps input points \(\mathbf{x} \in \mathbb{R}^d\) to the set of $2m$ sine waves:
\begin{equation}
    \gamma (\mathbf{x}) = 
    [\cos (2\pi \mathbf{B} \mathbf{x}),
    \sin (2\pi \mathbf{B} \mathbf{x})]^\mathrm{T},
\end{equation}
where each entry in \(\mathbf{B}\in \mathbb{R}^{m\times d}\) is sampled from \(\mathcal{N}(0, \sigma_{\text{FF}}^2)\). An example is shown in Test Case 2.

\subsection{Training algorithm}

We divide the model training into two phases, associated with \(\nnSDF\) and \(\nnPHY\) respectively, as described above. For the training of \(\nnSDF\), we initially perform some epochs employing the Adam optimizer, initialized with a learning rate of \(10^{-3}\). Subsequently, we transition the training process to full-batch training with the L-BFGS optimizer until the loss difference tolerance between two continuous epochs, which is set to be \(10^{-8}\), is achieved. The GeLU activation function is used in this work.

For the training of \(\nnPHY\), the combination of Adam and L-BFGS optimizer is still utilized. In this part, we consider hyperbolic tangent (tanh) activation functions.
We remark that the use of differentiable activation functions, such as tanh, maintain the framework's flexibility for integrating physics-informed loss terms.
The number of epochs considered in each test case is reported in Sec.~\ref{sec:results}.

\subsection{Evaluation metrics}

For the evaluation of \(\nnSDF\), we employ the Chamfer distance (CD). Each geometry in the test dataset is provided with a set of \(K\) data pairs \(\{\mathbf{x}_i, s_i\}_{i=1}^K\) for the learning of the shape code \(\hat{\mathbf{z}}_i\). With the learned shape code \(\hat{\mathbf{z}}_i\), a uniform mesh grid in the space \(\{\mathbf{x}_{i,j}\}_{i=1, j=1}^{N_x, N_y}\), where \(N_x\) and \(N_y\) are the discretization numbers along the x and y axes respectively, is passed through $\nnSDF$, which outputs the predicted SDF value on this mesh grid. The marching cubes method \cite{lorensen1998marching} is then used to identify the zero level set on the grid, and the resulting coordinate set \(S_1\) is considered as the reconstructed geometry by $\nnSDF$.
Given a point set \(S_2\), which contains boundary coordinates of reference geometries, the CD metric is simply the average of the nearest-neighbor distances for each point to the nearest point in the other point set with $N_{S_l}$ denoting the cardinality of \(S_l\) ($l = 1,2$):
\begin{equation}
    d_{CD}(S_1^i, S_2^i)=\frac{1}{N_{S_1^i}}\sum_{\mathbf{x}\in S_1^i}\min_{\mathbf{y}\in S_2^i}||\mathbf{x}-\mathbf{y}||^2_2 +
                     \frac{1}{N_{S_2^i}}\sum_{\mathbf{y}\in S_2^i}\min_{\mathbf{x}\in S_1^i}||\mathbf{x}-\mathbf{y}||^2_2,
\end{equation}
where $i=1,...,N_{t}$, and $N_t$ denotes the size of the test dataset. The average of CD is used for the evaluation of the test dataset which consists of $N_{t}$ samples, given as:
\begin{equation}
    \frac{1}{N_t}\sum_{i=1}^{N_{t}}d_{CD}(S_1^i, S_2^i).
\end{equation}
For the evaluation of \(\nnPHY\), we use the relative \(L_2\) error metric. This evaluation method is applied to a test dataset containing \(M_{t}\) samples, denoted as \(\{\Omega_i\}_{i=1}^{M_t}\). For each test sample, there are \(P_i\) test data points \(\{\mathbf{x}_i^j, \overline{s}_i^j\}_{j=1}^{P_i}\), which are randomly sampled within the geometry. These data points, along with the corresponding shape code \(\mathbf{z}_i\), are passed through \(\nnPHY\) to obtain the prediction. The relative \(L_2\) error for the entire test dataset is given by:

\begin{equation}
    \sqrt{
    \frac{\sum_{i=1}^{M_{t}}\sum_{j=1}^{P_i}
    |\nnPHY(\mathbf{x}_i^j,{\mathbf{z}}_i, \overline{s}_i^j; \paramsPHYhat)-
    \mathbf{u}(\mathbf{x}_i^j, \Omega_i)|^2}
    {\sum_{i=1}^{M_{t}}\sum_{j=1}^{P_i}
    |\mathbf{u}(\mathbf{x}_i^j, \Omega_i)|^2}}.
\end{equation}

\subsection{SDF-USM-Net workflow}
Using SDF-USM-Net, like other surrogate modeling approaches, comprises offline and online stages. 
Having introduced all the components, we now provide a summary in this section. Algorithm~\ref{alg:workflow-alg} presents the summarized workflow for both the offline and online stages.

\begin{algorithm}[H]
\caption{Workflow of SDF-USM-Net.}
\label{alg:workflow-alg}
\hspace*{\algorithmicindent} \textbf{Offline stage:\\}
\hspace*{\algorithmicindent} \textbf{Input:} 
the number of training shapes $N$ for $\nnSDF$, the number of training shapes $M$ for $\nnPHY$
\begin{algorithmic}[1]
\State Define the $N$ training shapes $\{\Omega_i\}_{i=1}^N$ for $\nnSDF$ 
\State Generate the training set $ G^{\mathrm{train}}:=\{G_i\}_{i=1}^N$
\State Train $\nnSDF$ by solving Eq. \eqref{eqn:SDF_training} until convergence
    \If{$M>N$}
        \For{$i = N+1$ to $M$}
            \State Define the training shape $\Omega_i$
            \State Compute the shape codes $\mathbf{z}_i$ by minimizing \eqref{eqn:shape_code_inference}   
        \EndFor

    \EndIf
\State Conduct FOM simulations for $M$ training shapes for $\nnPHY$ and generate the training set $\{U_i\}_{i=1}^M$
\State Launch the training of $\nnPHY$ based on the pre-trained set of shape codes until convergence in Eq. \eqref{eq:phy_training}
\end{algorithmic}
\hspace*{\algorithmicindent} \\
\hspace*{\algorithmicindent} \textbf{Online stage:\\}
\hspace*{\algorithmicindent} \textbf{Input:} 
$N_t$ new shapes\\
\hspace*{\algorithmicindent} \textbf{Output:} Inference of physical fields for $N_t$ new shapes
\begin{algorithmic}[1]
\State Convert the new set of shapes into data pairs in test dataset $G^{test}:=\{G_i\}_{i=1}^{N_t}$
\State Infer the shape codes $\mathbf{z}_i$ by minimizing \eqref{eqn:shape_code_inference}
\State Feed-forward the new shape codes \(\mathbf{z}_i\) together with query points \(\mathbf{x}\) and DF values through \(\nnPHY\)
\end{algorithmic}
\end{algorithm}

\subsection{Alternative shape encoding approaches} \label{sec:shape operator}

To demonstrate the effectiveness of the SDF-USM-Net for predicting physical quantities, we compare three types of shape-encoders \(P_g:\mathcal{G}\rightarrow \mathbb{R}^{k}\):
\begin{itemize}
    \item \textbf{Latent shape codes}. Shape codes learned by \(\nnSDF\), as described in Sec. \ref{sec:deepsdf}.
    \item \textbf{Explicit parameters}. In case an explicit parametrization of the elements of the space \(\mathcal{G}\) is available, we define \(P_g\) in such a way that the latent codes coincide with the geometrical parameters themselves. An example is provided in Test Case 2.
    \item \textbf{Geometrical Landmarks}. For geometries with complex shapes in which an explicit representation is not available, a straightforward choice is to take the coordinates of key points in the domain as landmarks, as proposed in \cite{regazzoni2022universal}. An example is provided in Test Case 1.
    
\end{itemize}

%% file: parts_data_preparation.tex
\section{Test cases}\label{sec:data}
In this Section, we introduce two test cases and offer an overview on the methodology employed for data generation.

\subsection{Test Case 1: Coronary bifurcation}
As a first test case (Test Case 1), we consider a test proposed in \cite{regazzoni2022universal}, consisting of the prediction of blood flow and pressure fields within a coronary bifurcation afflicted with stenosis. More precisely, our investigation focuses on a computational domain \(\Omega\) representing a two-dimensional cross-section of a coronary artery featuring a bifurcation. 

\subsubsection{Geometrical variability}
\begin{figure}
    \centering
    \includegraphics[width=0.7\linewidth]{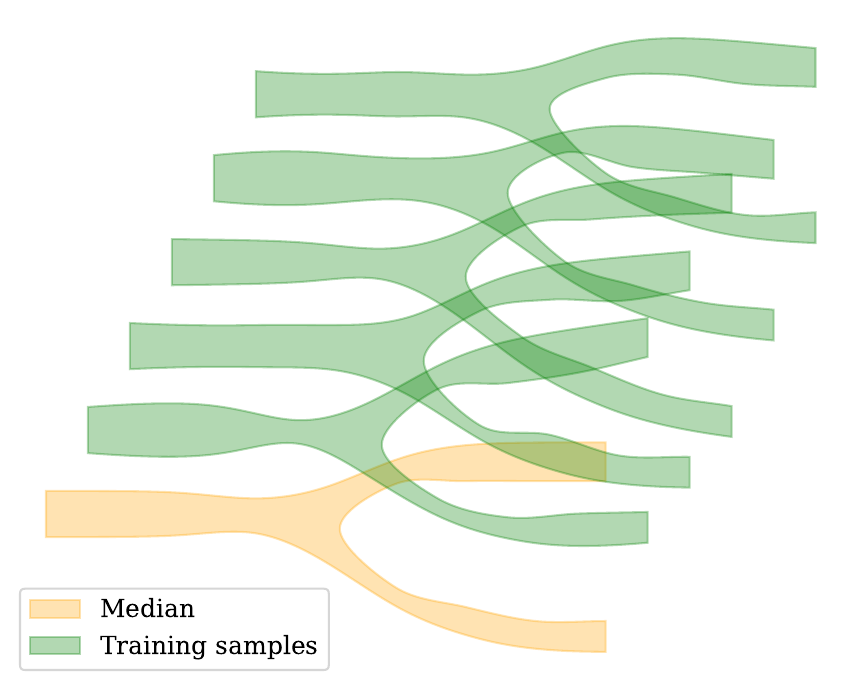}    
    \caption{Test Case 1: Representation of some of the geometries \(\Omega\in \mathcal{G}\) included in the training dataset in Test Case 1.}
    \label{fig:geo_var}
\end{figure}

We synthetically generate an extensive array of diverse computational domains, each representing virtual patients, by employing splines. These splines are derived by randomly varying their parameters within defined intervals, designed to emulate the realistic variability observed among patients. A subsample of the geometries obtained following this procedure is displayed in Fig~\ref{fig:geo_var}. All distances used in this test case are expressed in centimeters (cm). The length of all bifurcations along the \(x\)-axis is 4 cm.

\subsubsection{SDF calculation}
For the initial phase of training \(\nnSDF\), we construct the SDF data pairs \(G_i\) for each geometry. These samples comprise 2D points $\mathbf{x}_i^j$ paired with their corresponding SDF values $s_i^j$. 

\begin{figure}
    \centering
    \begin{subfigure}[b]{0.7\textwidth}
    \centering
    \includegraphics[width=\linewidth]{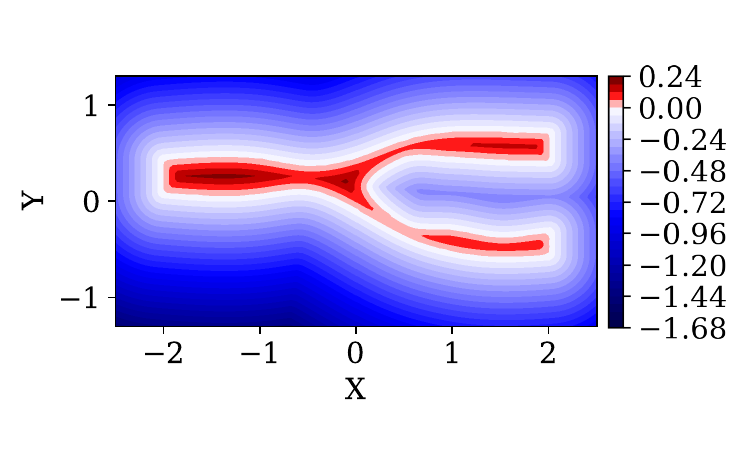}  
    \caption{}
    \end{subfigure}
    
    \vspace{0.3cm}

    \begin{subfigure}[b]{0.7\textwidth}
    \centering
    \includegraphics[width=\linewidth]{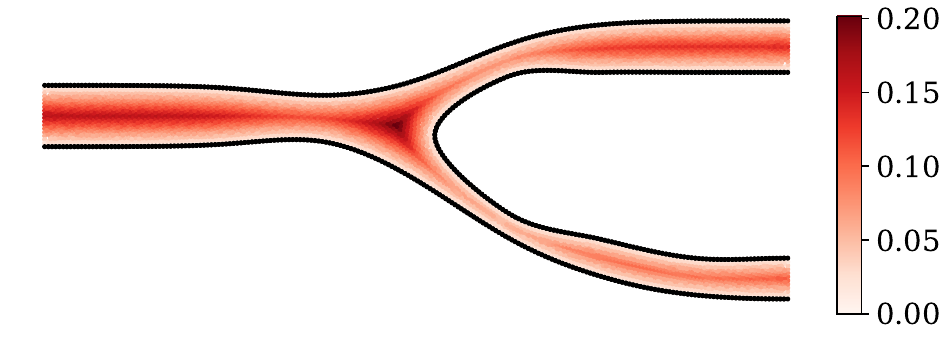}    
    \caption{}
    \end{subfigure}
    
    \caption{Test Case 1: (a) A sample of $SDF(\mathbf{x})$ on a uniform grid. (b) A sample of ${DF}(\mathbf{x})$ on an unstructured grid.}
    \label{fig:sdf}
\end{figure}

This process involves the following steps: Firstly, we randomly sample 8000 points along the geometry boundary. Next, we perturb these points with Gaussian noise: 4000 points with a variance of 0.033 and the remaining 4000 points with a variance of 0.0033. This results in a total of 8000 points representing the geometry. This sampling strategy aims to capture detailed SDF information near the boundary. Additionally, we uniformly sample points within an area slightly larger than the computational domain, with a discrete distance of 0.1, resulting in around 900 points. Subsequently, for each sampled point \(\mathbf{x}\), we determine its SDF value by identifying the closest point in the set of boundary points denoted as \(P\). This proximity search is facilitated using the KD-Tree algorithm\cite{10.1145/361002.361007}. Figure~\ref{fig:sdf}(a) shows a sample of generated SDF on a uniform grid.

For the second phase of training \(\nnPHY\), the DF values $\overline{s}$ are computed by determining the minimum distance from the query point $\mathbf{x}$ within the geometry to the physical boundary as shown in Fig.~\ref{fig:sdf}(b), which does not include inlet and outlets. Notice that this data generation process relies solely on the coordinates of the boundary. Consequently, the computation of the DF is applicable to both real and synthetic geometries.

\subsubsection{Geometrical landmarks}

\begin{figure}
    \centering
    \includegraphics[width=0.7\linewidth]{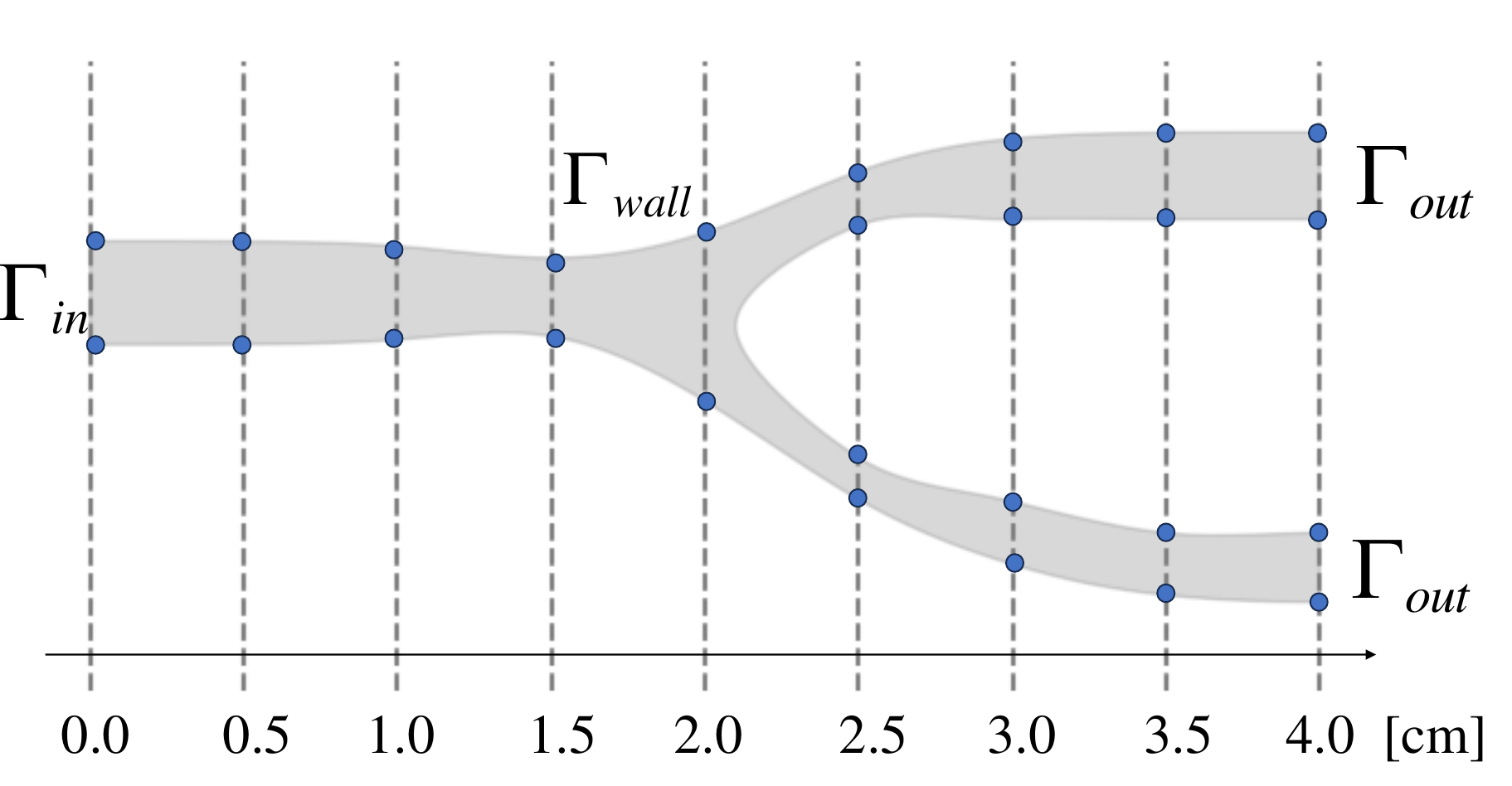}
    \caption{Test Case 1: geometrical landmarks.}
    \label{fig:landmarks}
\end{figure}

As a benchmark for the present test case, we consider the approach proposed in \cite{regazzoni2022universal}, wherein the domain shape variability is encoded through geometrical landmarks, due to the lack of an explicit parameterization of these geometries in practical problems. Specifically, we define landmarks as the coordinates $y$ of the vessel wall corresponding to a set of predefined coordinates $x$. In this test case, we will consider a set of landmarks containing 26 coordinates as shown in Fig~\ref{fig:landmarks}.

\subsubsection{Governing equation and solution generation}
 We consider the  following stationary Navier-Stokes model, describing the steady-state fluid flow in the coronary bifurcation
\begin{equation} \label{eqn:coronary}
    \left\{
        \begin{aligned}
        &- \nu \Delta \mathbf{v} +  ( \mathbf{v} \cdot \nabla) \mathbf{v} + \frac{1}{\rho} \nabla p = \mathbf{0}
        && \quad \text{in }\Omega,  \\
        & \mathbf{\nabla}  \mathbf{v} = 0
        && \quad \text{in } \Omega,  \\
        & \mathbf{v} = \mathbf{v}_{\text{in}}
        && \quad \text{on } \Gamma_{\text{in}}, \\
        & \mathbf{v} = \mathbf{0}
        && \quad \text{on } \Gamma_{\text{wall}}, \\
        & \nu \frac{\partial \mathbf{v}}{\partial \mathbf{n}} - \frac{1}{\rho} p \, \mathbf{n} = \mathbf{0}
        && \quad \text{on } \Gamma_{\text{out}}, \\
        \end{aligned}
    \right.
\end{equation}
where \(\nu = 4.72\ \text{mL}^2\text{s}^{-1}\) is the kinetic viscosity of blood and \(\rho=1060\ \text{kg\ m}^{-3}\) is the density. At the inlet, we consider a parabolic profile with a maximum velocity equal to $14$ cm s\(^{-1}\). 

\begin{figure}
    \centering
    \includegraphics[width=0.8\linewidth]{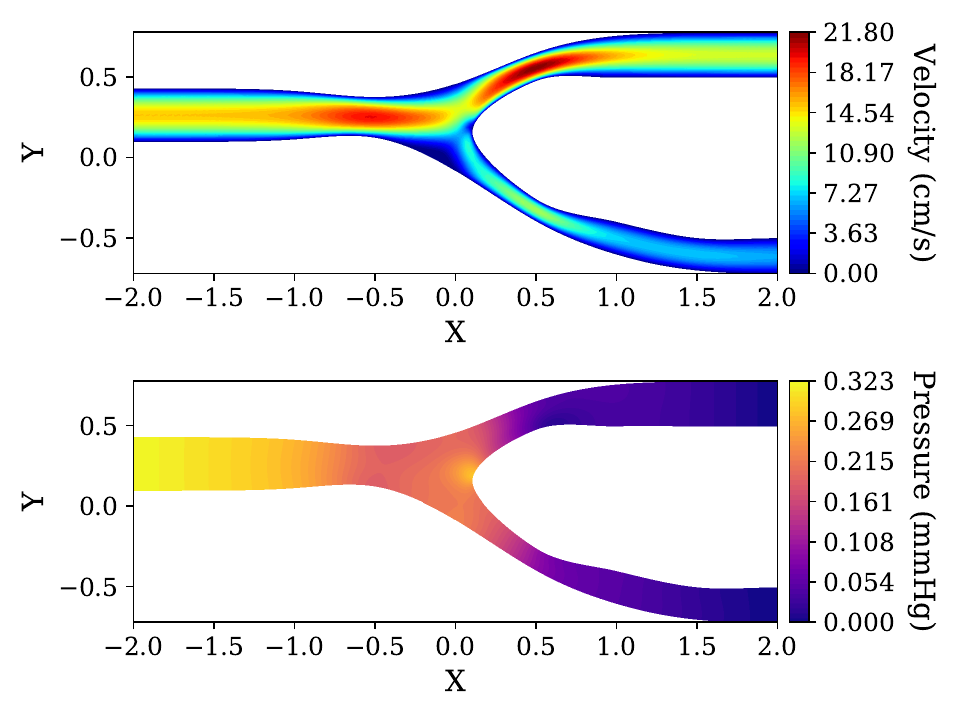}
    \caption{Test Case 1: An example of solution.}
    \label{fig:phy_sol}
\end{figure}

We solve Eq.~\eqref{eqn:coronary} by the finite element method implemented in the free/open-source software FEniCS \cite{Baratta_DOLFINx_the_next_2023}. To generate training data in the second phase, we utilize a Taylor-Hood Finite Element solver. For spatial discretization, we adopt triangular computational meshes with a spatial resolution of approximately 
$h=0.2 \ \text{mm}$. An example solution is shown in Fig.~\ref{fig:phy_sol}.   

\subsection{Test Case 2: Elastic plate}
Test Case 2 involves predicting the displacement of an elastic body with a square shape, featuring geometrical and topological changes defined by random configurations of holes positioned at arbitrary locations. 

\subsubsection{Geometrical variability}
We synthetically generate a wide range of computational domains, each representing plates with one or two randomly positioned holes. This test scenario aims to assess the effectiveness of our proposed model in addressing challenges related to the topological variation of the main features of a geometry. 

Furthermore, in this test case, we compare the accuracy achieved with the proposed shape encoding technique to the accuracy obtained with an explicit exact parametrization, to assess how closely the proposed method approaches this best-case scenario.
Hence, for any given geometry in the set \(\mathcal{G}\), the vector \(\mathbf{\mu}_g=[x_1, y_1, R_1,...,x_i, y_i, R_i, N]\), where $x_i$ and $y_i$ are the coordinates of the circle center, $R_i$ is the radius of the circle, and \(N\) is the number of the holes, are used as explicit parameters. We define the operator \(P_g(\Omega):=\mathbf{\mu}_g\) to represent the geometries in the given dataset. 

Here we present two sub-cases:
\begin{itemize}
    \item \textbf{Test Case 2a: fixed topology.} To tackle the challenges posed by discrete positional changes in geometry shapes, we train our model using a dataset exclusively comprising plates with two randomly positioned holes of varying radii, parameterized by \(\mathbf{\mu}_g=[x_1, y_1, R_1,x_2, y_2, R_2, 2]\), where the topology (i.e., the number of holes) remains fixed.
    \item \textbf{Test Case 2b: topology changes.} To address the challenges presented by topology changes, we train our model using a dataset consisting of plates with either one or two randomly positioned holes. In this case, the parameterization \(\mathbf{\mu}_g=[x_1, y_1, R_1,...,x_i, y_i, R_i, N]\), where \(N\in \{1,2\}\), accounts for variations in both the positions and the number of holes, reflecting the changes in topology.
\end{itemize}

\subsubsection{SDF calculation}
For the initial phase of training our SDF-USM-Net, since each plate shares the same square frame shape, we exclude the square boundary during the generation of SDF values. This exclusion mitigates the influence of the uniform square feature on the model's ability to discern variations among circular shapes.

For the second phase of training our SDF-USM-Net, the $DF$ values are calculated as the minimum distance between a given point inside the physical domain and all the physical boundaries.

\subsubsection{Governing equation and solution generation}
A square plate made of an isotropic and homogeneous hyperelastic material under finite deformation is considered. We denote by \(\zeta\) the mapping from the reference configuration to the deformed one, so that the reference coordinate $\boldsymbol{X}$ and the deformed one $\boldsymbol{x}$ are linked by:
\begin{equation}    \boldsymbol{x}=\zeta(\boldsymbol{X})=\boldsymbol{X}+{\boldsymbol{u}(\boldsymbol{X})},
\end{equation}
where $\boldsymbol{u}\colon \Omega \to \mathbb{R}^2$ denotes the displacement field.
We consider traction forces on the left and right boundaries ($\Gamma_\text{right}$ and $\Gamma_\text{left}$), and normal spring boundary conditions on the top and bottom boundaries ($\Gamma_\text{top}$ and $\Gamma_\text{bottom}$).
Hence, the equilibrium configuration is given by the solution of the following PDE, where \(\mathbf{P}\) is the first Piola-Kirchhoff stress and \(\mathbf{N}\) represents the outward normal unit vector in the initial configuration:
\begin{equation}
    \left\{
        \begin{aligned}
        &\nabla \cdot \mathbf{P}=\mathbf{0}
        && \text{in }\Omega,  \\
        & \mathbf{P} \mathbf{N} ={t}\mathbf{N}
        && \text{on } \Gamma_\text{right} \cup \Gamma_\text{left},  \\
        & \mathbf{P} \mathbf{N} + K \mathbf{N} \otimes \mathbf{N} \boldsymbol{u} = \mathbf{0}
        && \text{on } \Gamma_\text{top} \cup \Gamma_\text{bottom},  \\
        \end{aligned}
    \right.
\end{equation}
with traction force $t = 0.1$, and spring stiffness $K = 0.01$.
The constitutive law for the hyperelastic material is written as:
\begin{equation}
    \begin{aligned}
        &\mathbf{P}=\frac{\partial \psi (\mathbf{F})}{\partial \mathbf{F}}, \qquad
        &\mathbf{F}=\nabla\zeta = I + \nabla\boldsymbol{u},
    \end{aligned}
\end{equation}
where \(\bf{F}\) presents the deformation gradient, and \(\psi\)  denotes the strain energy density of a certain material. In this case, we consider a Neo-Hooken material response, and the strain energy density is given by:
\begin{equation}
    \psi (\mathbf{F})
    =
    \frac{1}{2}\mu \,(I_1-2)
    -\mu \, \log (J)
    + \frac{1}{2}\lambda \, (\log(J))^2
    ,
\end{equation}
where the first principal invariant is determined by \(I_1=\text{trace}(\mathbf{C})\), the right Cauchy-Green tensor \(\mathbf{C}\) is defined as \(\mathbf{C}=\mathbf{F}^T\mathbf{F}\), \(J=\text{det}(\mathbf{F})\) is the determinant of the deformation gradient. The Lam\'e parameters are set to \(\mu=0.5\) and \(\lambda=1.3\). 

\begin{figure}
    \centering
    \includegraphics[width=0.75\linewidth]{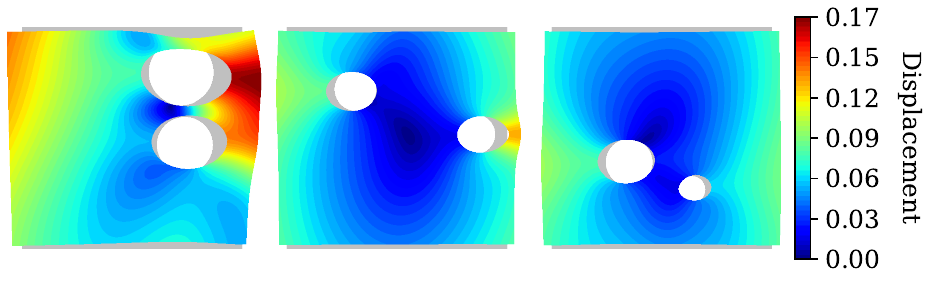}
    \caption{Test Case 2: Three examples of the circle locations and the corresponding displacement of the elastic plates.}
    \label{fig:sol_2}
\end{figure}

To generate training data, we employ a linear Finite Element solver based on Fenics \cite{Baratta_DOLFINx_the_next_2023}.  For space discretization, we consider triangular computational meshes with a space resolution of nearly \(h=0.3\). Figure \ref{fig:sol_2} shows examples of displacement solutions contained in the training set.

%% file: parts_results.tex
\section{Results} \label{sec:results}
In this Section, we report our numerical investigation of the performances of the proposed SDF-USM-Net on the two test cases outlined in Sec. \ref{sec:data}.

\subsection{Test Case 1: Coronary bifurcation}
In this test case, our training set comprises \(N=500\) distinct geometries.

\subsubsection{\(\nnSDF\) training results} \label{sec:recons_1}
 We employ an FCNN with four inner layers, each comprising 32 neurons, chosen based on a thorough hyperparameter tuning procedure (see \ref{sec:arch-test} for more details). 
 The smoothed truncation function is utilized to let the training of \(\nnSDF\) focus on the details of the geometries, and the truncation distance \(\beta\) is set to $0.1$. The regularization parameter \(\sigma\) is set to \(10^{2}\).

The learning rate for the shape codes is set to \(10^{-5} \times B\), where \(B\) represents the number of shapes in one batch. For each geometry, we select \(K_i=8900\), for $i=1,...,N$, randomly generated points for training \(\nnSDF\). During training with the Adam Optimizer, we randomly subsample 625 SDF samples (out of 8900 available points) for each shape in a batch, over 1000 epochs. Subsequently, we transition to full-batch training with the L-BFGS Optimizer, utilizing 1000 subsamples for each shape until the loss difference reaches the tolerance of \(10^{-8}\). 

\begin{figure}
    \centering
    \includegraphics[width=0.4\linewidth]{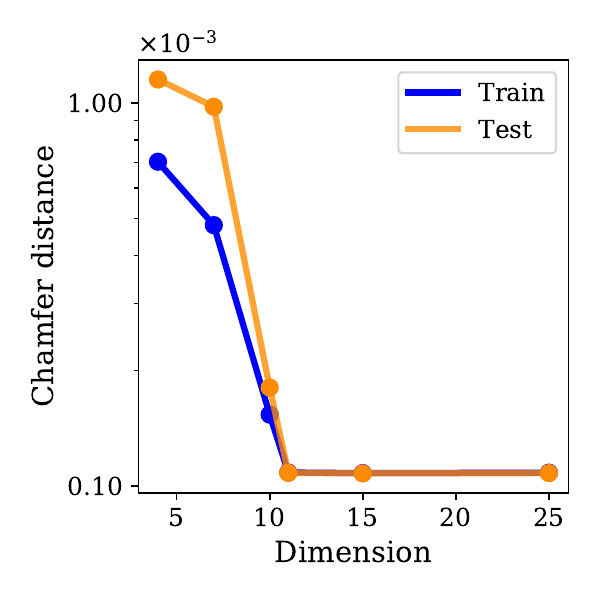}
    \caption{Test Case 1: impact of latent dimensions. Using Chamfer distance as an evaluation metric for the reconstructed shapes, illustration of the relationship between Chamfer distance and latent dimensions, ranging from 4 to 25.}
    \label{fig:cd_dim}
\end{figure}

The dimension of the latent space $k$ emerges as a critical hyperparameter highly impacting the learning capability of the model. In this case, due to the dataset comprising variations in 11 parameters, we conducted experiments around \(k=11\). Figure \ref{fig:cd_dim} illustrates the Chamfer distances obtained by varying the latent dimensions from 4 to 25. As observed, an increase in dimension correlates with a decrease in Chamfer distance. However, when the dimension reaches or exceeds 11, the errors on both the training and test sets become nearly equal reaching down to the magnitude of \(1\times 10^{-4}\), with no further reduction in error observed with increasing dimension. 

\begin{figure}
    \centering
    \includegraphics[width=0.85\linewidth]{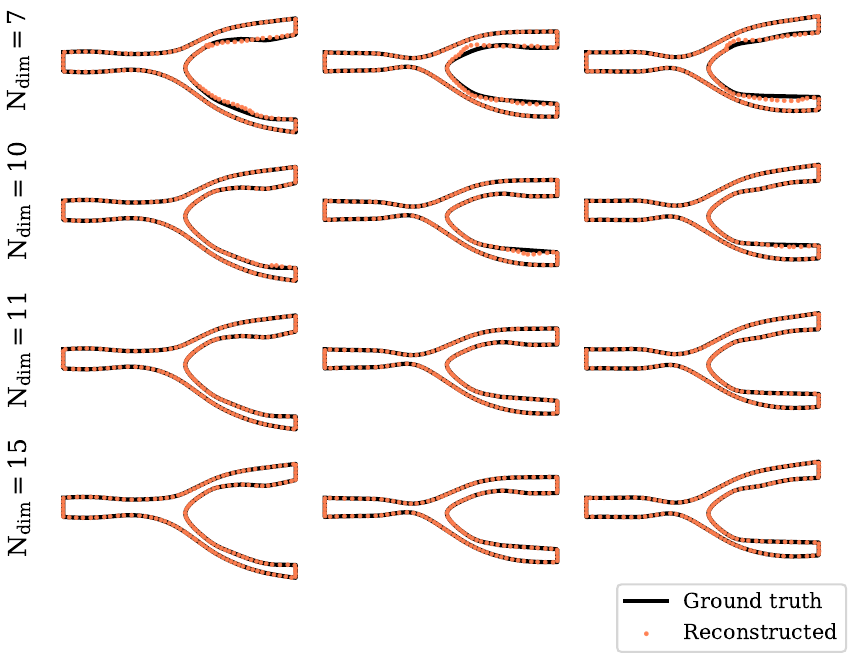}
    \caption{Test Case 1: impact of latent dimension. With each column representing an example from the test dataset, the results represent the reconstructed geometries with an increasing number of latent dimensions (7, 10, 11, 15).}
    \label{fig:recons_1}
\end{figure}

Figure~\ref{fig:recons_1} presents reconstructions for three samples, showcasing dimensions of 7, 10, 11, and 15, respectively. It is evident that, for dimensions 7 and 10, visible differences persist between the reconstructed shapes and the ground truth, with the case of dimension 7 exhibiting a more dispersed spatial distribution. Conversely, when the dimension is 11 or higher, the reconstructed shapes align closely with the true shapes, indicating that an 11-dimensional space is sufficient to effectively encode the shape distribution in the latent space.
Since these geometries have been synthetically generated by varying 11 parameters, we can conclude that, through the analysis of Figure \ref{fig:recons_1}, we have been able to re-discover the intrinsic dimension of the set of shapes $\mathcal{G}$ in an automated manner.

\subsubsection{Fluid velocity and pressure field prediction}
For the training of \(\nnPHY\), we utilize \(M_i=1000\) randomly generated points for each computational domain. An FCNN with four inner layers is employed, comprising 20, 15, 10, and 5 neurons, respectively. 
To train the FCNN weights and biases, we employ the Adam optimizer for 200 iterations and then switch to the L-BFGS optimizer until convergence with the loss difference tolerance set to \(10^{-8}\). 

\begin{figure}
    \centering
    \begin{subfigure}[b]{0.78\textwidth}
        \centering
        \includegraphics[width=\linewidth]{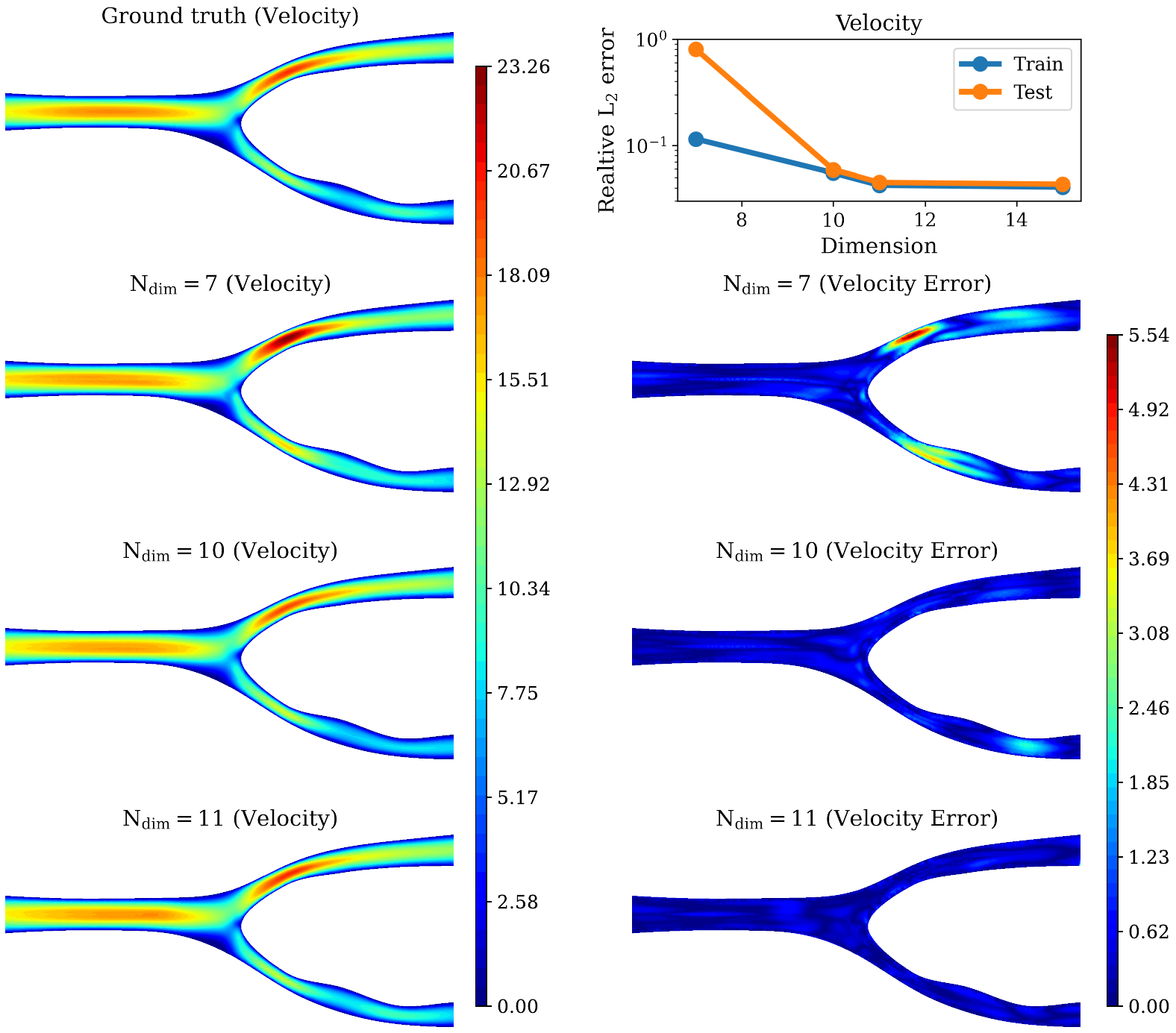}
        \caption{}
    \end{subfigure}
    \\
    \begin{subfigure}[b]{0.78\textwidth}
        \centering
        \includegraphics[width=\linewidth]{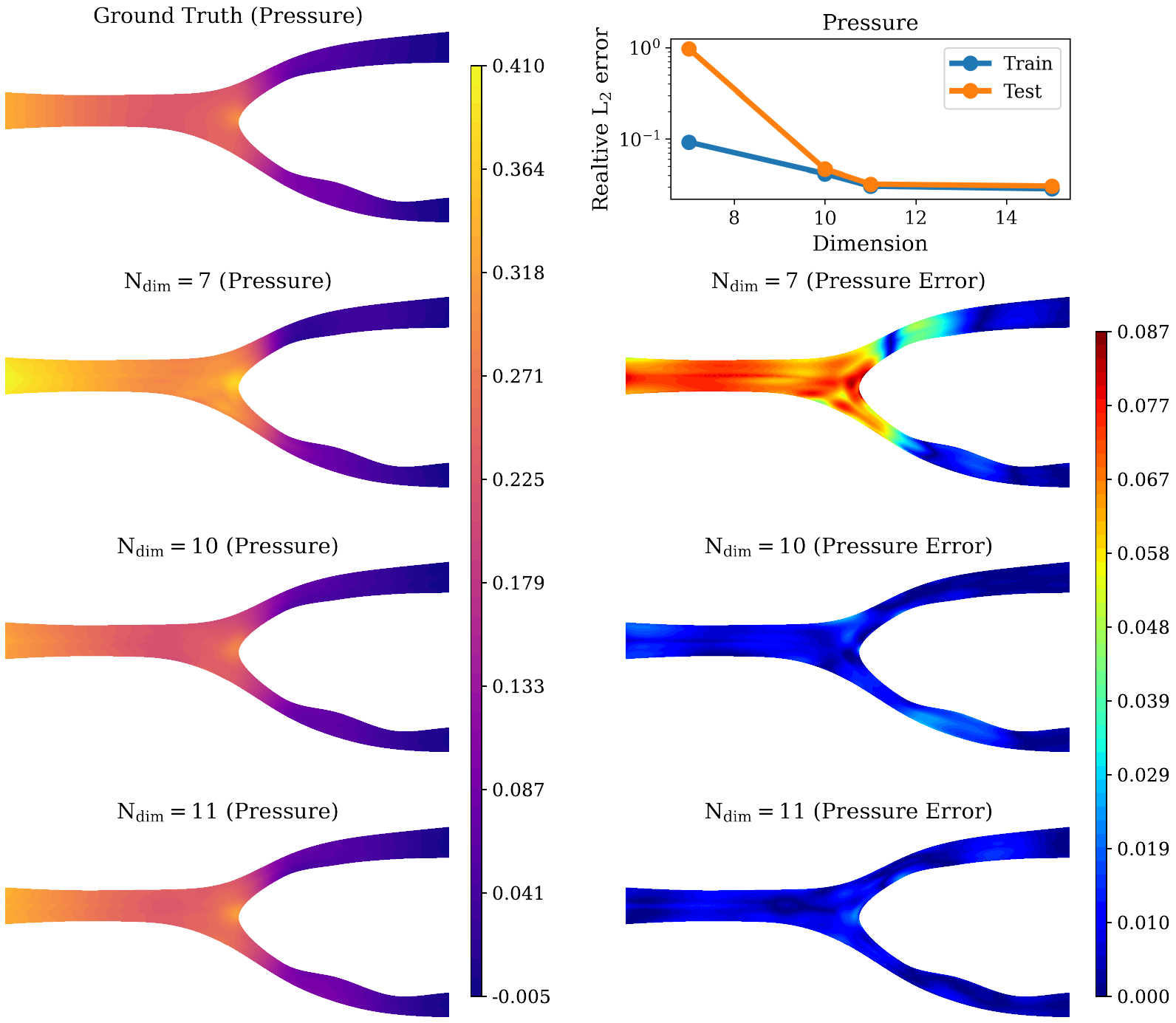}   
        \caption{}
    \end{subfigure}
    
    \caption{Test Case 1: impact of latent space dimensions. We compare, for a sample belonging to the test dataset, the results obtained by using SDF-USM-Net with an increasing number of latent states (7, 10, 11). For (a) velocity magnitude, and (b) pressure, the top-left corner reports the FOM solution, followed by the pressure prediction. The right column reports the absolute error concerning the FOM solution. The relative \(L_2\) error of Test Case 1 as a function of the number of latent dimensions is shown at the top-right corner.}
    \label{fig:surro_1_dim}
\end{figure}

\paragraph{Impact of latent space dimension}
Based on the experiments of shape reconstruction discussed in Section \ref{sec:recons_1}, it is evident that the dimension of the latent space plays a crucial role in training the SDF-USM-Net. Here, we focus on dimensions around 11 and select 7, 10, 11, and 15.
To illustrate the impact of the latent dimension on the reconstruction of the flow field, we randomly selected a sample from the test dataset for visualization in Fig. \ref{fig:surro_1_dim}. It is apparent from the visualization that as the dimension increases, the prediction error, particularly that associated with the pressure, exhibits a significant decrease. The relative \(L_2\) error as a function of latent dimension showcases a trend similar to the Chamfer distance. Notably, there is no significant decrease in error observed with dimensions exceeding $11$. Based on these experiments, we choose the latent space dimension of $11$ for the following experiments.

\paragraph{Impact of input features}
To highlight the effectiveness of the proposed SDF-USM-Net, we consider four types of input features:
\begin{itemize}
    \item \textbf{USM-Net with landmarks}, which receives as input two spatial coordinates, and the geometrical landmarks.
    \item \textbf{SDF-USM-Net with shape codes}, which receives as input two spatial coordinates, and the latent expression of shape codes. 
    \item \textbf{USM-Net with landmarks and $DF$ values}, which receives as input two spatial coordinates, the geometrical landmarks, and the $DF$ values.     
    \item \textbf{SDF-USM-Net with shape codes and $DF$ values}, which receives as input two spatial coordinates, the latent expression of shape codes, and the $DF$ values. 
\end{itemize}

\begin{figure}
    \centering
    \begin{subfigure}[b]{0.49\linewidth}
    \centering
         \includegraphics[width=\linewidth]{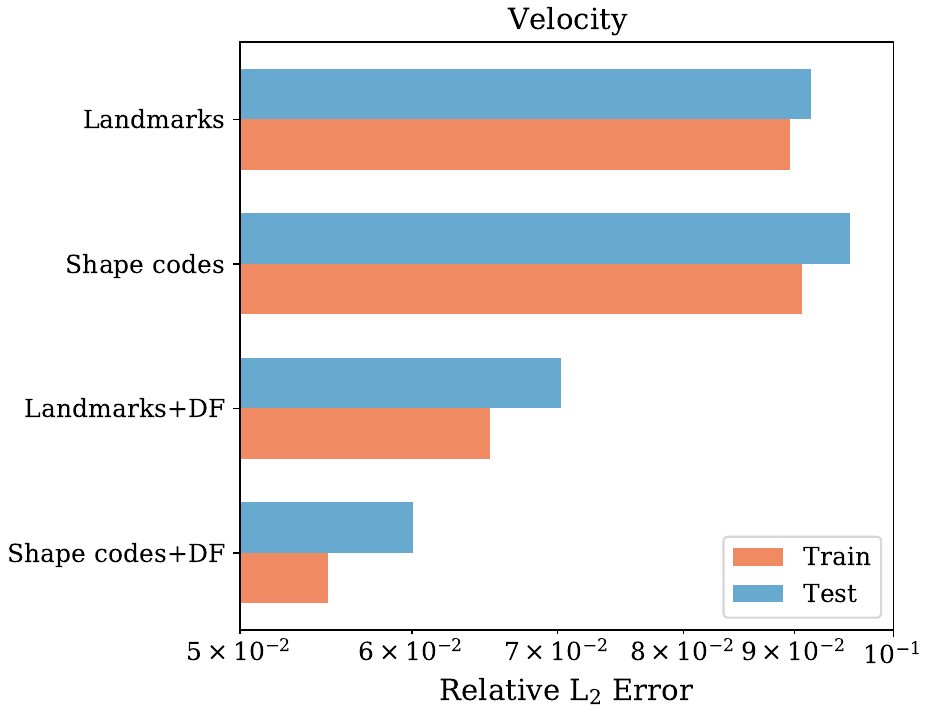}
         \caption{}
         \label{fig:in_v}
    \end{subfigure}
    \begin{subfigure}[b]{0.49\linewidth}
    \centering
         \includegraphics[width=\linewidth]{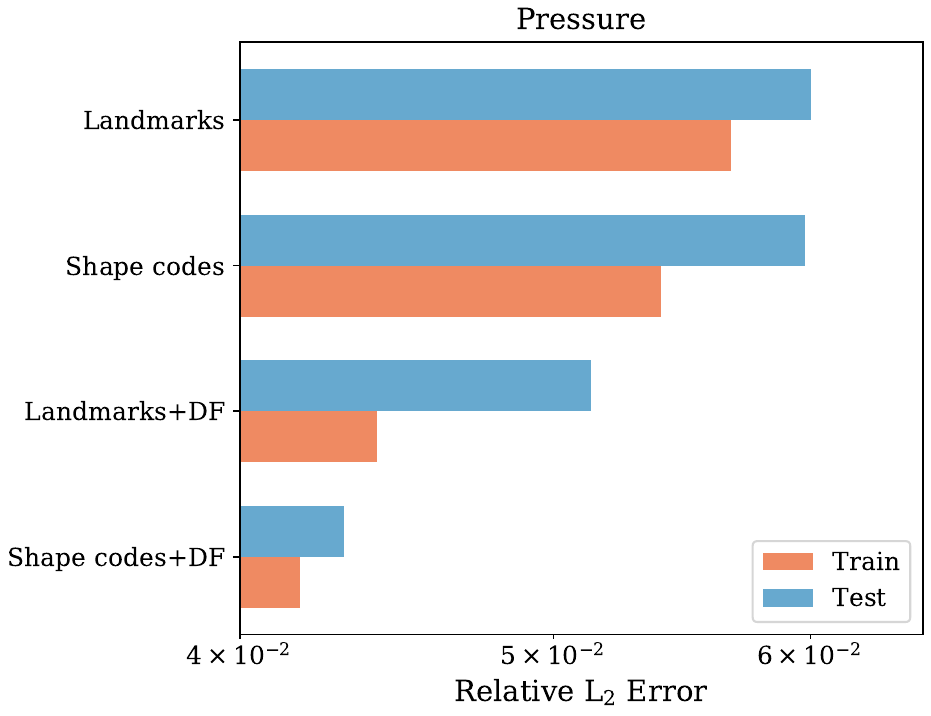}
         \caption{}
         \label{fig:in_p}
    \end{subfigure}
    \caption{Test Case 1: bar plots of the errors on the test dataset obtained with geometric landmarks and shape codes, and with or without $DF$ as additional feature. Left: error on the velocity field; right: error on the pressure field.}
    \label{fig:in_type}
\end{figure}

We present bar plots depicting the relative \(L_2\) error for velocity and pressure in Fig. \ref{fig:in_type}. From these results, a surrogate model trained with shape codes exhibits comparable error levels to the model trained with landmarks, which encompass all variable parameters of the geometry. Notably, shape codes can outperform landmarks when training is supplemented with the $DF$ values.
Furthermore, regardless of whether landmarks or shape codes are used as input, training combined with $DF$ values consistently enhances model performance.

\begin{figure}
    \centering
    \includegraphics[width=1.0\linewidth]{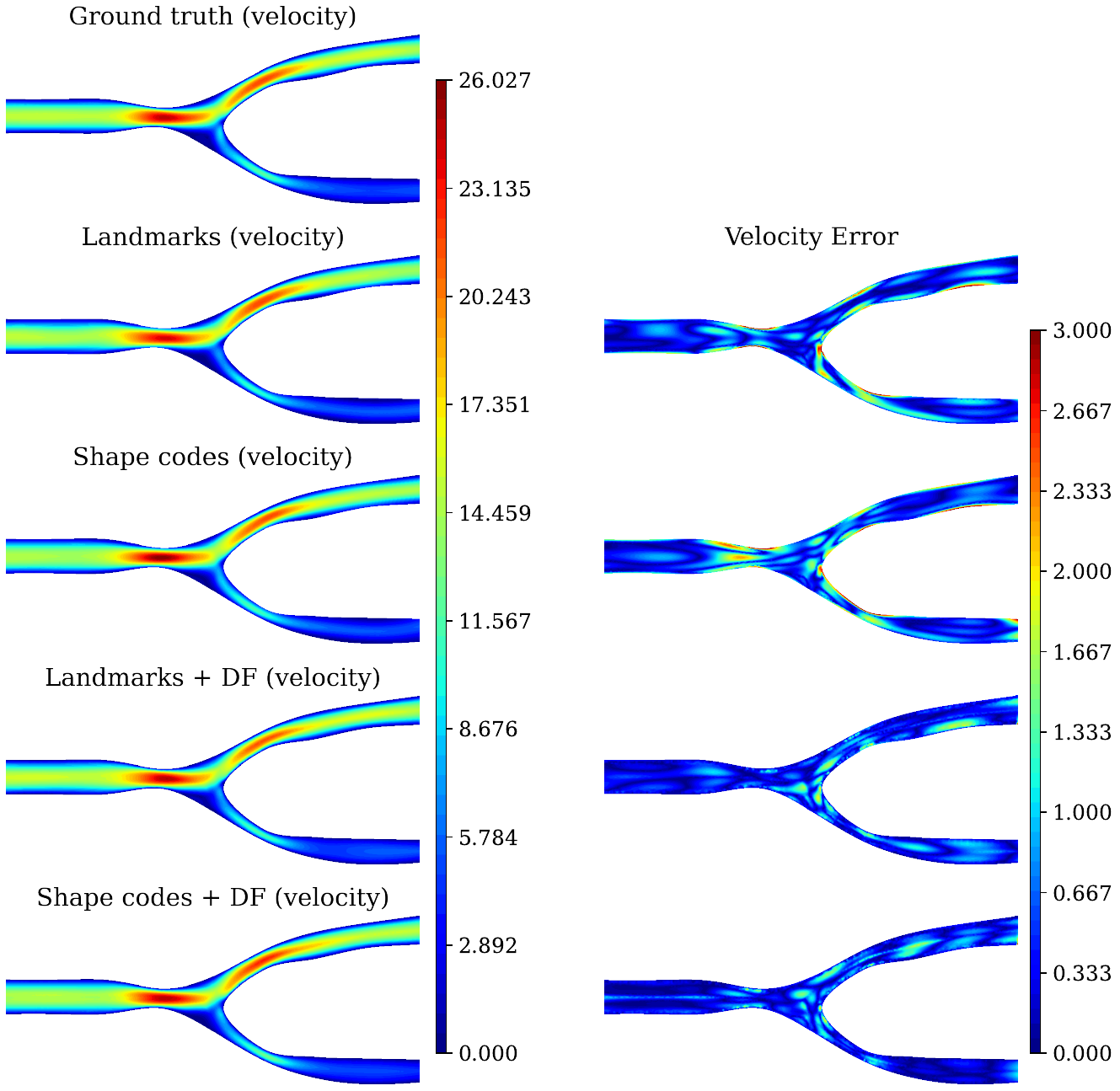}
    \caption{Test Case 1: impact of input features (velocity). We compare, for a sample belonging to the test dataset, the results obtained by using the original USM-Net (with geometrical landmarks) and the proposed SDF-USM-Net (with shape codes), with and without $DF$ values. The top-left corner reports the FOM solution for the velocity magnitude, followed by the prediction obtained with different input features. The right column reports the absolute error associated with the FOM solution.}
    \label{fig:surro_1_v}
\end{figure}

To gain further insight into the influence of different input features, we compare the predicted velocity magnitude obtained with all the considered architectures for a representative geometry belonging to the test dataset in Fig. \ref{fig:surro_1_v}. The prediction results indicate that the error without $DF$ values is notably larger near the boundary, particularly in regions adjacent to the stenosis, where significant shape variations occur. Conversely, models training with $DF$ values exhibit more precise predictions near the boundary.

\begin{figure}
    \centering
    \includegraphics[width=1.0\linewidth]{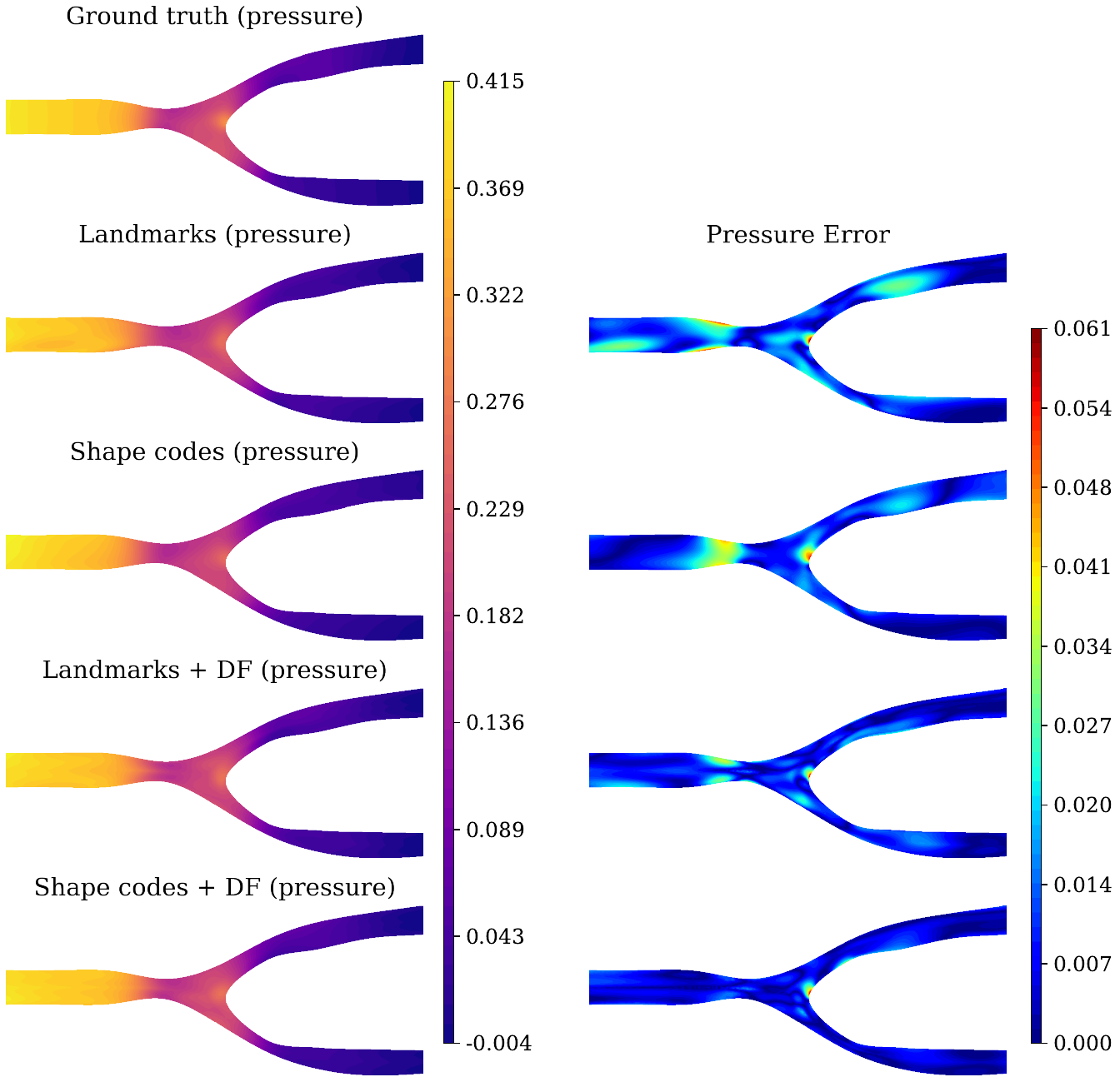}
    \caption{Test Case 1: impact of input features (pressure). We compare, for a sample belonging to the test dataset, the results obtained by using the original USM-Net (with geometrical landmarks) and the proposed SDF-USM-Net (with shape codes), with and without $DF$ values. The top-left corner reports the FOM solution for the pressure, followed by the prediction obtained with different input features. The right column reports the absolute error associated with the FOM solution.}
    \label{fig:surro_1_p}
\end{figure}

The comparative results for the same example on pressure are presented in Fig. \ref{fig:surro_1_p}. The disparity between different models is more pronounced here. Although all models can learn at least the coarse distribution of pressure, the SDF-USM-Net trained with shape codes and $DF$ values demonstrates superior prediction capabilities.

This test case highlights the superior performance in predicting the fluid flow solution when providing $DF$ values as additional features, especially close to the boundary of the domain.

\subsection{Test Case 2: Elastic plates}

In this section, we first report the results for Test Case 2a (fixed topology) and then for Test Case 2b (variable topology).

\subsubsection{Test Case 2a: fixed topology}
In this test case, our training set comprises \(N=1000\) distinct geometries. 
For \(\nnSDF\), we employ an FCNN with four inner layers, each containing 64 neurons. 
In this test case, we do not make use of a truncation function in training $\nnSDF$. This choice is explained later, and the impact of the truncation function is investigated.

The shape codes \(\mathbf{z}_i\) are initialized from the same normal distribution as defined in Section \ref{sec:deepsdf}. The dimension of the latent space is set to 5. We incorporate Fourier feature mapping and centralization, as described in Sections \ref{sec:mapping} and \ref{sec: central}, respectively, in the training of \(\nnSDF\). All other training details of \(\nnSDF\) remain consistent with those described in Section \ref{sec:recons_1}.

\begin{figure}
    \centering
    \includegraphics[width=1\linewidth]{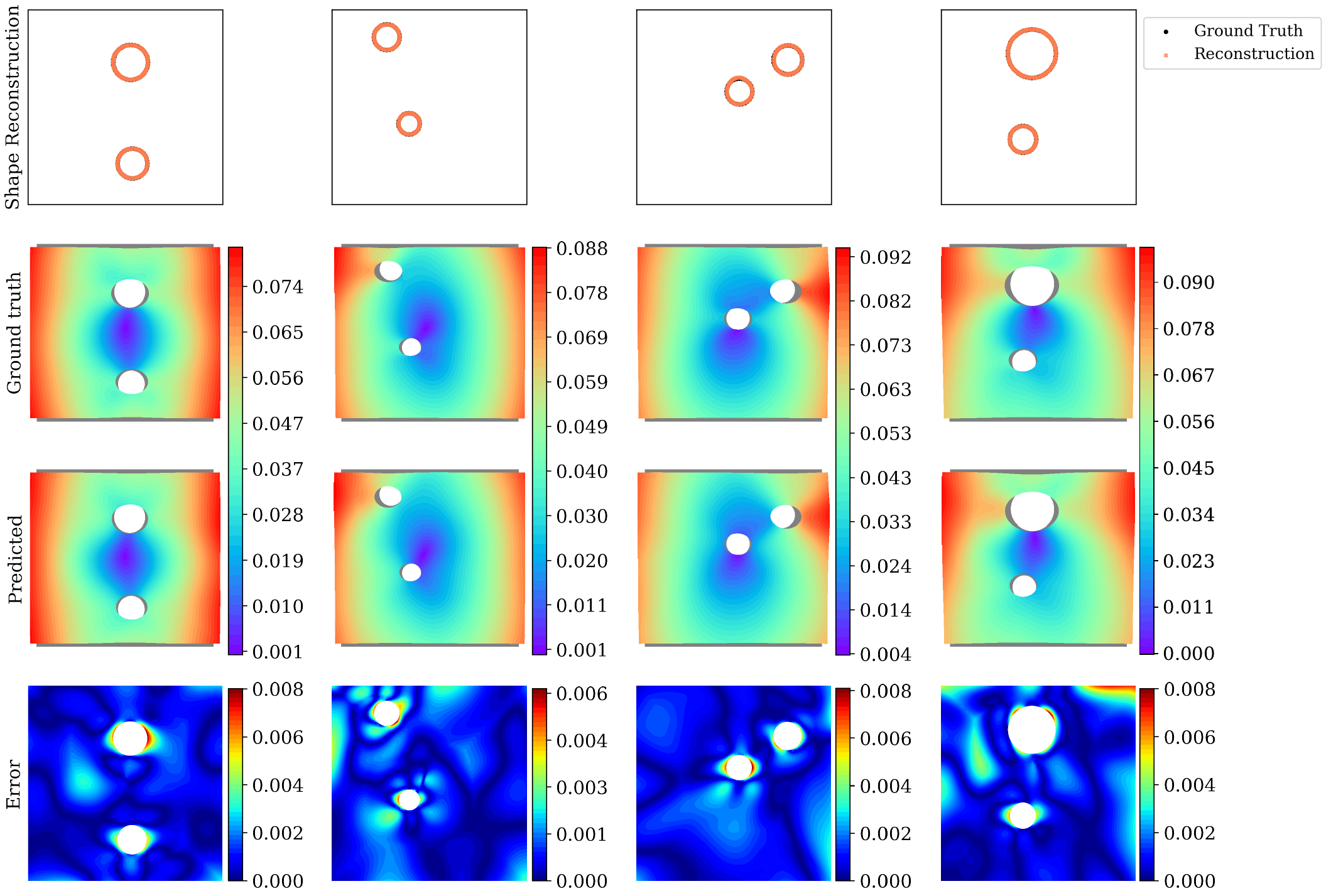}
    \caption{Test Case 2a: comparison of SDF-USM-Net reconstruction on four numerical solutions from the test set: from the best case scenario (first column, with a relative \(L_2\) error of 0.0256) to the worst case scenario (last column, with a relative \(L_2\) error of 0.1330).}
    \label{fig:example_2}
\end{figure}

For the training of \(\nnPHY\), we utilize \(M_i=1000\) randomly generated points for each computational domain. An FCNN with four inner layers is employed, comprising 30, 25, 20, and 15 neurons, respectively.

\(\nnSDF\) was trained on a dataset containing 750 samples and tested on 500 samples without experiencing overfitting (Chamfer distance on training dataset: \(8.513\times 10^{-5}\), on test dataset: \(8.515\times 10^{-5}\). These combined 1250 shapes were then divided into 80\% training and 20\% test data for training \(\nnPHY\). Some examples of contours of displacement resulting from SDF-USM-Net compared with the FOM results are shown in Fig. \ref{fig:example_2}. SDF-USM-Net can effectively learn and predict continuous changes in displacement over a given range of geometrical configurations.

\begin{figure}
    \centering
    \begin{subfigure}[b]{0.7\textwidth}
        \centering
        \includegraphics[width=\textwidth]{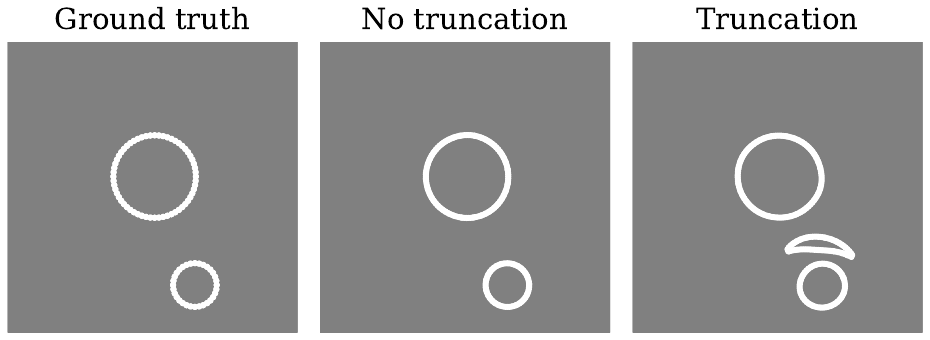}
        
    \end{subfigure}
    \vspace{0.3cm}
    \begin{subfigure}[b]{0.7\textwidth}
        \centering
        \includegraphics[width=\textwidth]{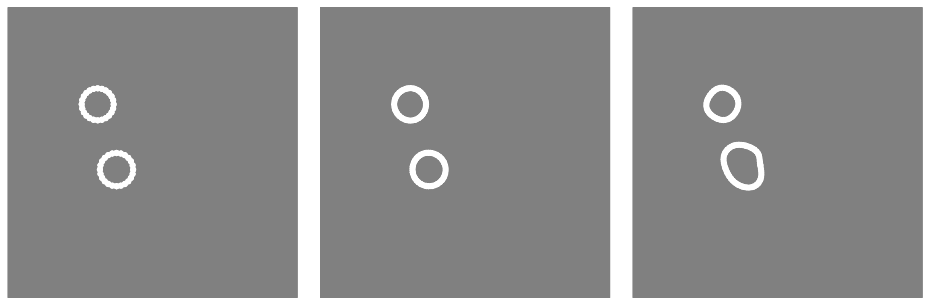}
    \end{subfigure}

    \caption{Test Case 2a: impact of truncation function. Two examples of corrupted reconstruction by applying truncation function in \(\nnSDF\).}
    \label{fig:clamp}
\end{figure}

\paragraph{Truncation function}
Whether to use the truncation function for the training of \(\nnSDF\) has a significant impact on the model's ability to capture macroscopic changes in shapes. Figure \ref{fig:clamp} illustrates failed reconstructions when the truncation function \(\text{clamp}_\beta (\cdot)\) is applied. Each reconstruction, involving the training of shape codes \(\mathbf{z}_i\) in the inference stage, is repeated 5 times to select the best reconstruction from different initializations. The results suggest that using the truncation function leads to a higher likelihood of getting trapped in local minima during the training of shape codes in the latent space which results in the failed reconstruction of the geometry.

\begin{figure}
    \centering
    \includegraphics[width=0.7\linewidth]{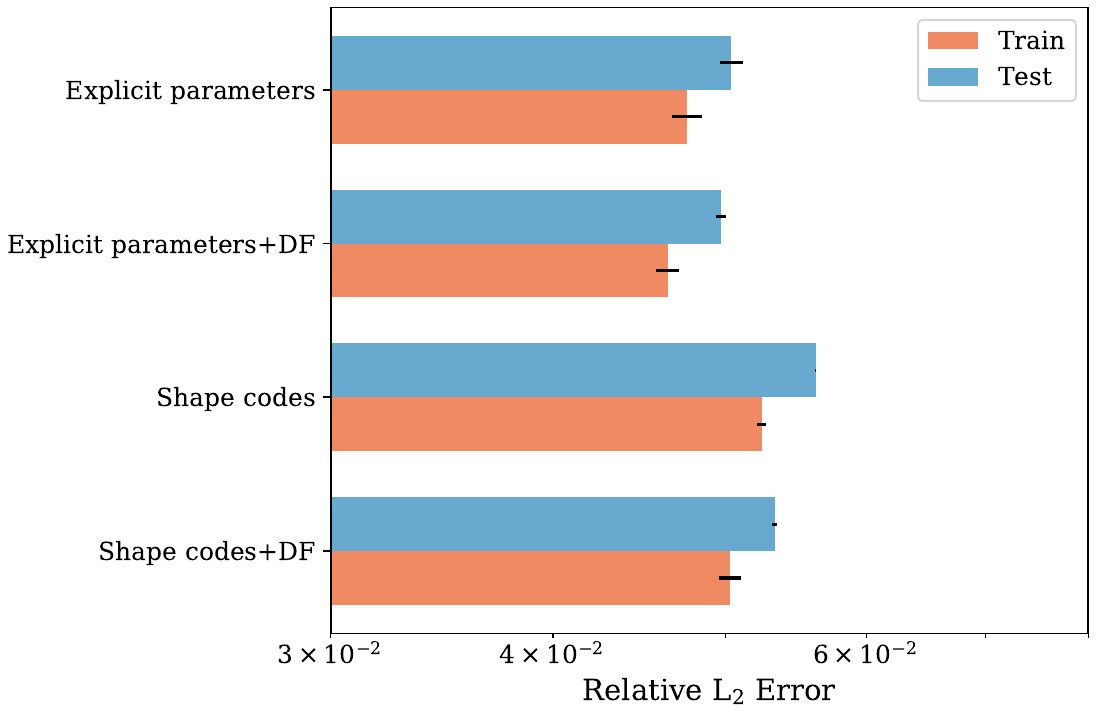}
    \caption{Test Case 2a: barplots of the errors on the training and test dataset obtained with different input features. The barplots refer to 3 training runs obtained starting from different random initializations of ANN weights and biases.}
    \label{fig:in_type_dx}
\end{figure}

\paragraph{Impact of input features}
We benchmark our proposed method with results obtained by replacing shape codes with the exact parameters used to generate the geometries. This idealized situation provides a way to obtain a reference level, a maximum level of accuracy achievable in principle, since explicit knowledge of the parameters allows the geometry to be encoded exactly.
The training results are summarized in Fig \ref{fig:in_type_dx}.

All experiments were repeated three times, and the mean and variance were calculated. These results show that training with input structure that consists of shape codes, remarkably, can achieve almost the same accuracy as training with the explicit knowledge of the parameters that were used to generate the geometry at hand. Moreover, training with input structure that consists of $DF$ values demonstrates slight improvements compared to training without $DF$ values.

\begin{figure}
    \centering
    \begin{subfigure}[b]{\textwidth}
        \centering
        \includegraphics[width=\linewidth]{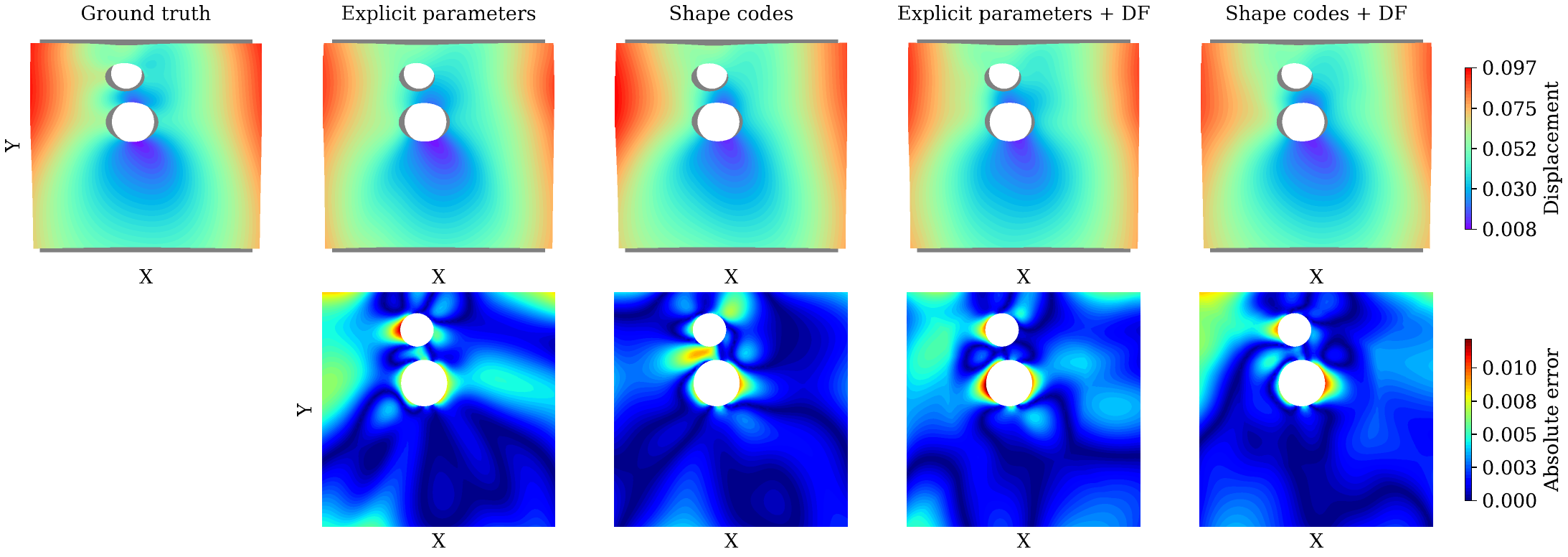}
    \end{subfigure}
    \\
    \begin{subfigure}[b]{\textwidth}
        \centering
        \includegraphics[width=\linewidth]{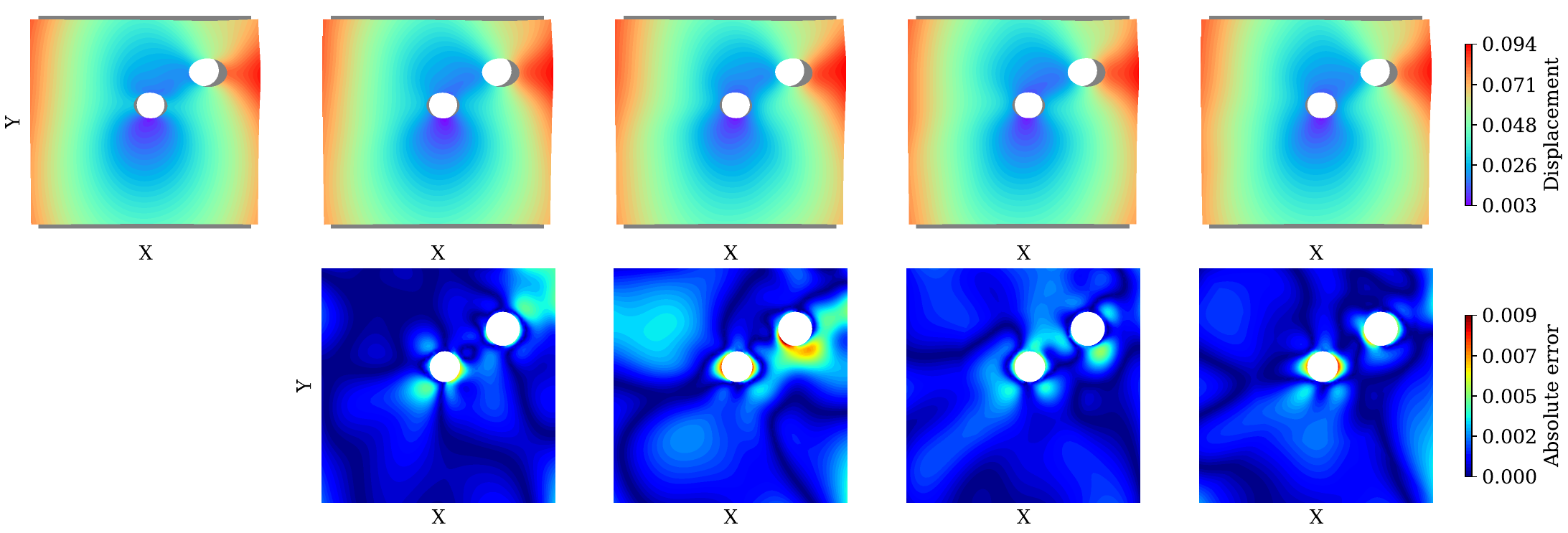}   
    \end{subfigure}
    \\
    \begin{subfigure}[b]{\textwidth}
        \centering
        \includegraphics[width=\linewidth]{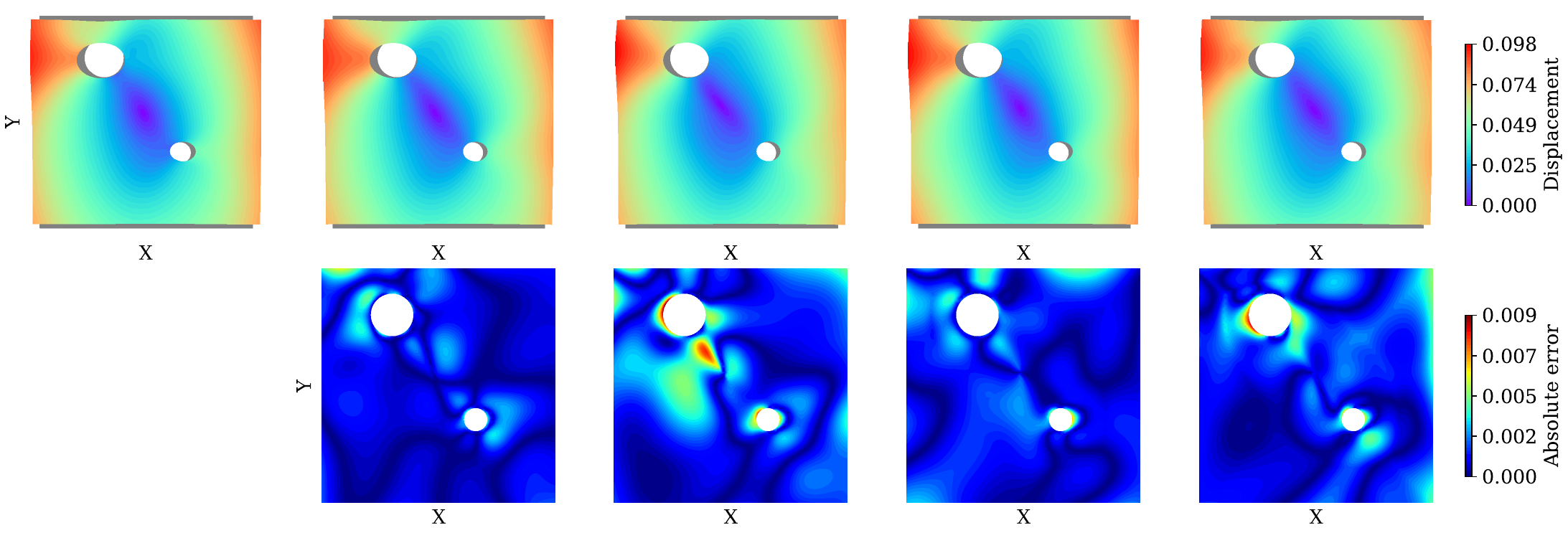}   
    \end{subfigure}
    \caption{Test Case 2a: comparison of predictions from different input structures. From top to bottom displays three samples in the test dataset.}
    \label{fig:plot_intype}
\end{figure}

Additionally, the statistical results from multiple experiments suggest that training with shape codes offers slightly greater stability to the random initialization of trainable parameters compared to training with explicit parameters on the test dataset. In Fig. \ref{fig:plot_intype} we show the displacement fields predicted by different input structures.

\subsubsection{Test Case 2b: variable topology}
\begin{figure}
    \centering
    \includegraphics[width=\linewidth]{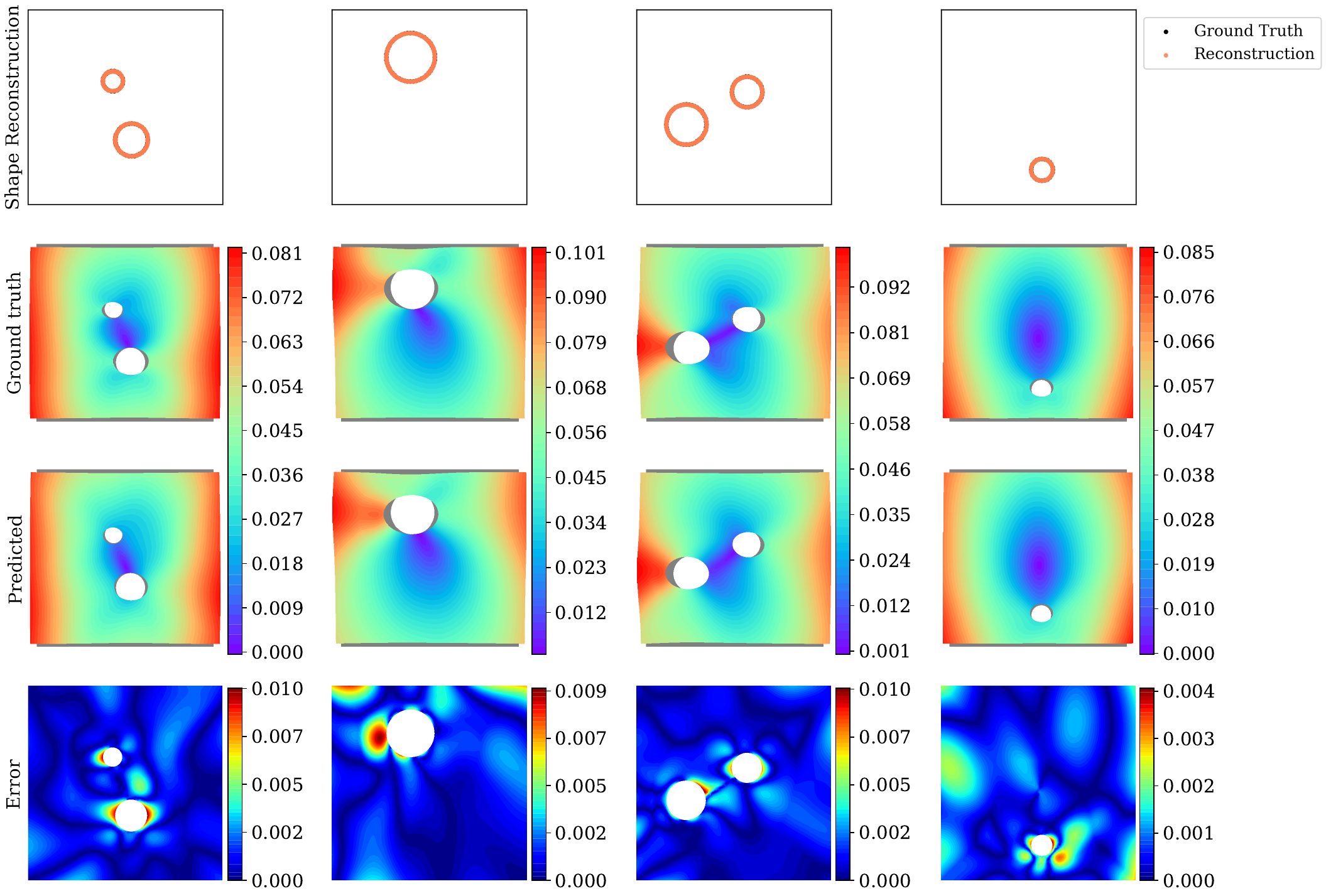}
    \caption{Test Case 2b: comparison of SDF-USM-Net reconstruction on four numerical solutions from the test set.}
    \label{fig:topology_change}
\end{figure}

We now consider the case of topological changes. The training setting is similar to Test Case 2a, while the latent dimension is changed to 8 due to the additional complexity related to the number of holes on the elastic plate. Figure \ref{fig:topology_change} shows examples of the reconstructed shape of the holes and predicted displacement. The reconstructions of displacement align consistently with the ground truth.

\begin{figure}
    \centering
    \includegraphics[width=0.83\linewidth]{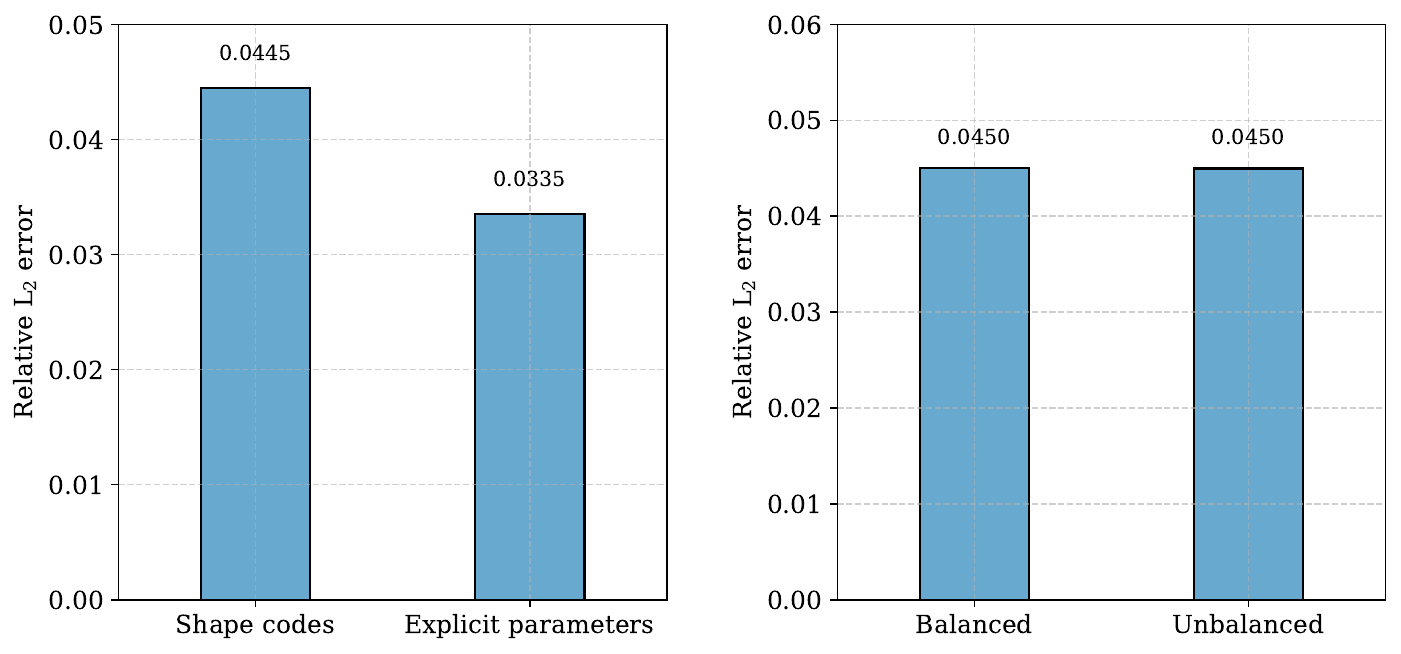}
    \caption{Test Case 2b: (a) Relative \(L_2\) error of different input structures. (b) Relative \(L_2\) error of balanced and unbalanced dataset.}
    \label{fig:in_balance}
\end{figure}

\paragraph{Impact of input features}
We further compare the relative \(L_2\) error in displacement between using shape codes and explicit parameters in Fig. \ref{fig:in_balance}a. The results indicate that the model's predictions can achieve a similar level of accuracy to that observed in Test Case 2a even in the presence of topology changes.

\paragraph{Impact of unbalanced dataset}
We evaluate SDF-USM-Net on both balanced and unbalanced datasets. In the balanced dataset, 50\% of the samples contain one hole, while the remaining 50\% contain two holes. Conversely, in the unbalanced dataset, approximately 30\% of the samples contain two holes, and the remaining 70\% contain one hole. The model is evaluated on the same test set, which contains balanced data.

From the results presented in Fig. \ref{fig:in_balance}b, it is observed that the model achieves almost the same accuracy for both datasets which suggests that our model also demonstrates robustness to the unbalanced distribution of data.

%% file: parts_discussion.tex
\section{Discussion} \label{sec:discussion}

We introduced a new method, termed SDF-USM-Net and built upon USM-Net \cite{regazzoni2022universal}, which integrates DeepSDF to learn the solution manifold of PDEs in a data-driven and geometry-universal manner.

A first NN, termed \(\nnSDF\), automatically learns a compact representation of shapes in a latent space, without the need for explicit encoder construction. The model merely requires pairs of coordinates and corresponding SDF values to characterize a given shape in a low-dimensional latent space, simplifying database construction (and reducing its dimensionality) compared to methods like auto-encoders. By implicitly representing geometries in the Euclidean space, the model is capable of extracting shape codes for geometries not included in the training set, demonstrating very good generalization capabilities. This renders our model an exceptionally flexible technique capable of addressing a broad spectrum of real-world applications, encompassing both shape and topology changes.

Unlike standard approaches that reconstruct a high-dimensional discretization of the output, such as evaluating the solution on the nodes of a computational mesh or of a Cartesian grid, our approach is meshless in nature. The reconstruction NN, termed \(\nnPHY\), is queried independently for each point in space. The meshless nature of SDF-USM-Nets, coupled with automatic feature extraction in low-dimensional latent spaces, makes it possible to learn the solution map underlying PDE models with lightweight architectures and with low training costs.

Specifically, our approach demonstrates the potential for clinical applications, as illustrated in Test Case 1. Given the complexity of the geometry, artificially defined geometrical landmarks may lack accuracy in representing it. In contrast, we notice that shape codes accurately encode the geometrical features of the domain, enhancing representational accuracy. With an 11-dimensional latent space, in fact, the model achieves a prediction error for the flow field of 6.01\%, showing an improvement over predictions based on more than 20 geometrical landmarks, where the error is 7.03\%. Additionally, we have presented tests of our model on geometry databases undergoing topological changes, as shown in Test Case 2. With the inclusion of an additional input layer based on Fourier features, the model achieves a better representation capacity, achieving a Chamfer distance of $8.69\times 10^{-5}$. 
Our results indicate that the proposed method, even when relying on learned shape codes rather than the knowledge of explicit geometric parameters, can achieve a level of accuracy nearly equivalent to the ideal scenario with known parameters, with errors of 5.33\% and 4.97\%, respectively.
Considering that an explicit parametrization may be not available in many practical problems, SDF-USM-Net exhibits a remarkable performance.

A current limitation of this work is that it does not consider the case of time-dependent inputs and variable physical parameters, which will be the subject of future research. However, we emphasize that the proposed method demonstrates flexibility, automation, and accuracy within its current scope.

%% file: parts_appendix.tex
\appendix
\section{Comparison study}
In this section, we list the detailed experiment setups and parameter searching for test cases in Sec. \ref{sec:results}.

\subsection{Test Case 1: coronary bifurcation}
We explored the NN architecture design space and optimization strategy through the following experiments.
\paragraph{Network depth and width}\label{sec:arch-test}
We tested network depths ranging from 3 to 5 layers and widths from 32 to 64 neurons for training $\nnSDF$. The results of these experiments, conducted on 300 training samples, are presented in Table \ref{table:archi}. In the following experiments for Test Case 1, the depth is set to 4 with width set to 32.

\begin{table}[ht]
\centering
\begin{tabular}{|c|c|}
\hline
\textbf{Depth$\times$width} & \textbf{Chamfer distance} \\ \hline
3$\times$32 & 1.075$\times 10^{-3}$ \\ \hline
4$\times$32 & 3.117$\times 10^{-4}$ \\ \hline
5$\times$32 & 6.298$\times 10^{-4}$ \\ \hline
4$\times$48 & 8.014$\times 10^{-4}$ \\ \hline
4$\times$64 & 6.307$\times 10^{-4}$ \\ \hline

\end{tabular}
\caption{Network architecture vs Chamfer distance on the test dataset.}
\label{table:archi}
\end{table}

\paragraph{Optimization strategy}
With the network depth and width determined, we trained \(\nnSDF\) using both a one-step optimization strategy (Adam only) and a two-step optimization strategy (Adam followed by L-BFGS), as shown in Table \ref{table:optim}. The NN in this experiment was trained on 500 samples. Our findings indicate that the two-step optimization strategy achieves higher accuracy compared to the one-step optimization strategy.

\begin{table}[ht]
\centering
\begin{tabular}{|c|c|c|}
\hline
\textbf{$k$} & \textbf{Adam} & \textbf{Adam+L-BFGS} \\ \hline
15 & 2.479$\times 10^{-4}$ & 1.389$\times 10^{-4}$ \\ \hline
25 & 1.743$\times 10^{-4}$ & 1.388$\times 10^{-4}$ \\ \hline
\end{tabular}
\caption{Chamfer distance on the test dataset after one-step vs two-step optimization.}
\label{table:optim}
\end{table}

\paragraph{Activation function}
We tested several types of activation functions, and the results on the test dataset are summarized in Table \ref{table:activation_functions}. Our findings indicate that \(\nnSDF\) using the GeLU activation function outperforms other activation functions in the current application scenario.

\begin{table}[ht]
\centering
\begin{tabular}{|c|c|c|c|}
\hline
\textbf{Activation function} & \textbf{ReLU} & \textbf{Tanh} & \textbf{GeLU} \\ \hline
 Chamfer distance & 1.389$\times 10^{-4}$ & 2.144$\times 10^{-4}$ & 1.091$\times 10^{-4}$ \\ \hline
\end{tabular}
\caption{Performance comparison of activation functions}
\label{table:activation_functions}
\end{table}

\subsection{Test Case 2: elastic plate}
The discrete and topological changes in the geometry present significant challenges to the training of $\nnSDF$. In this section, we evaluate the effects of centralization and Fourier feature mapping on improving the model's performance.
\paragraph{Centralization}
The centralization method, as defined in Sec. \ref{sec: central}, is designed to improve the model's robustness to geometric variations by aligning the centroid of each shape. This alignment ensures that the network can focus on the intrinsic properties of the shapes rather than their arbitrary positions. The results on the test dataset presented in Table \ref{table:central} demonstrate a notable improvement in accuracy when the centralization method is applied, indicating its effectiveness in handling the discrete changes in geometry.

\begin{table}[ht]
\centering
\begin{tabular}{|c|c|c|}
\hline
\textbf{Coordinate input} & \textbf{No centralization} & \textbf{Centralization}\\ \hline
 Chamfer distance & 1.430$\times 10^{-4}$ & 8.763$\times 10^{-5}$ \\ \hline
\end{tabular}
\caption{Performance improvement of centralization}
\label{table:central}
\end{table}

\paragraph{Fourier feature mapping}
Fourier feature mapping, described in Sec. \ref{sec:mapping}, enhances the network's ability to capture low-frequency variations in the input data.  The results on the test dataset shown in Table \ref{table:fourier} illustrate the performance gains achieved through this technique.

\begin{table}[ht]
\centering
\begin{tabular}{|c|c|c|c|c|}
\hline
\textbf{Input mapping} & \textbf{No mapping} & \textbf{m=16} & \textbf{m=8} & \textbf{m=4}\\ \hline
 Chamfer distance & 7.44$\times 10^{-4}$ & 2.33$\times 10^{-4}$ & 8.72$\times 10^{-5}$ & 8.69$\times 10^{-5}$\\ \hline
\end{tabular}
\caption{Performance improvement of Fourier feature mapping}
\label{table:fourier}
\end{table}

%% file: parts_acknowledgements.tex
F.R. and S.P. has received support from the project PRIN2022, MUR, Italy, 2023--2025, P2022N5ZNP ``SIDDMs: shape-informed data-driven models for parametrized PDEs, with application to computational cardiology''. 
The present research is part of the activities of ``Dipartimento di Eccellenza 2023--2027'', MUR, Italy, Dipartimento di Matematica, Politecnico di Milano. 
F.R. and S.P. acknowledge their membership to INdAM GNCS - Gruppo Nazionale per il Calcolo Scientifico (National Group for Scientific Computing, Italy). 
J. Z. and L. Z. has received support from the National Natural Science Foundation of China (grant No. 92371102).